\documentclass[10pt]{article}

\usepackage[nocompress]{cite}
\usepackage{geometry}                		
\usepackage{graphicx}	
\usepackage[table]{xcolor}				
\usepackage{amssymb}
\usepackage{amsmath}
\usepackage{mathtools}
\usepackage{amsthm}
\usepackage{multirow}
\usepackage[hidelinks]{hyperref}
\usepackage{bm}
\newcommand{\cca}{\cellcolor{blue!10}}
\newcommand{\ccb}{\cellcolor{red!10}}

\usepackage{amsthm}
\usepackage{algorithmic}
\usepackage{algorithm}
\usepackage{booktabs}
\usepackage{lineno}
\usepackage{julia-listings-format}
\lstdefinelanguage{JuliaLocal}{
    language = Julia, 
    morekeywords = [3]{toeplitz, construct_shifted_Tncm, construct_Bn, construct_M_P, construct_positive_real_eigenvalues, construct_all_eigenvalues,construct_B}, 
    morekeywords = [2]{LinearAlgebra, GenericSchur}, 
}

\makeatletter
\newcommand{\doubletilde}[1]{{%
  \mathpalette\double@tilde{#1}%
}}
\newcommand{\double@tilde}[2]{%
  \sbox\z@{$\m@th#1\tilde{#2}$}%
  \ht\z@=.9\ht\z@
  \tilde{\box\z@}%
}
\makeatother
\usepackage{amsthm}
\newtheorem{rmrk}{Remark}

\newtheorem{conj}{Conjecture}

\begin{document}

\title{On the eigenvalues of Toeplitz matrices with two off-diagonals}

\author{Sven-Erik Ekström\thanks{Department of Information Technology, Uppsala University}\\{\small\href{mailto:sven-erik.ekstrom@it.uu.se}{sven-erik.ekstrom@it.uu.se}} \and David Meadon\footnotemark[1]\\{\small\href{mailto:david.meadon@it.uu.se}{david.meadon@it.uu.se}}}

\date{}

\maketitle

\begin{abstract}
  \noindent
Consider the Toeplitz matrix \(T_n(f)\) generated by the symbol \(f(\theta)=\hat{f}_r e^{\mathbf{i}r\theta}+\hat{f}_0+\hat{f}_{-s} e^{-\mathbf{i}s\theta}\), where \(\hat{f}_r, \hat{f}_0, \hat{f}_{-s} \in \mathbb{C}\) and \(0<r<n,~0<s<n\). For \(r=s=1\) we have the classical tridiagonal Toeplitz matrices, for which the eigenvalues and eigenvectors are known. Similarly, the eigendecompositions are known for \(1<r=s\), when the generated matrices are ``symmetrically sparse tridiagonal''.

In the current paper we study the eigenvalues of \(T_n(f)\) for \(1\leq r<s\), which are ``non-symmetrically sparse tridiagonal''.
We propose an algorithm which constructs one or two ad hoc matrices smaller than \(T_n(f)\), whose eigenvalues are sufficient for determining 
the full spectrum of \(T_n(f)\). The algorithm is explained through use of a conjecture for which
examples and numerical experiments are reported for supporting it and for clarifying the presentation. Open problems are briefly discussed.
\end{abstract}

{\small \textbf{\textit{Keywords:}} Toeplitz matrix sequences, eigenvalues, tridiagonal

}

\section{Introduction}
\label{sec:int}
We say that a function \(f\in L^1(-\pi,\pi)\) generates a Toeplitz matrix \(T_n(f)\in\mathbb{C}^{n\times n}\) if
\begin{align*}
T_n(f)=[\hat{f}_{i-j}]_{i,j=1}^{n}=
\left[
\begin{array}{ccccccc}
\hat{f}_0&\hat{f}_{-1}&\hat{f}_{-2}&\dots&\hat{f}_{1-n}\\
\hat{f}_1&\hat{f}_0&\hat{f}_{-1}&\ddots&\vdots\\
\hat{f}_2&\ddots&\ddots&\ddots&\vdots\\
\vdots&\ddots&\ddots&\ddots&\hat{f}_{-1}\\
\hat{f}_{n-1}&\dots&\dots&\hat{f}_{1}&\hat{f}_0
\end{array}
\right],
\end{align*}
where \(\hat{f}_k\) are the Fourier coefficients of \(f\),
\begin{align*}
\hat{f}_{k}=\frac{1}{2\pi}\int_{-\pi}^{\pi}f(\theta)e^{-\mathbf{i}k\theta}\mathrm{d}\theta,\qquad \mathbf{i}^2=-1,\quad k\in\mathbb{Z}.
\end{align*}
The function \(f\) is called the generating symbol and is also a symbol for the singular value distribution in the Weyl sense, according to classical results by Szeg\H{o}~\cite{Szeg__1920}, Avram~\cite{Avram_1988}, Parter~\cite{Parter_1986}, B\"{o}ttcher~\cite{Bottcher2012-lf}, Silbermann~\cite{Bottcher2012-lf}, Tilli~\cite{Tilli_1998}, Tyrtyshnikov~\cite{Tyrtyshnikov_1998},  Serra-Capizzano~\cite{Capizzano_2003,Serra_Capizzano_2006}. Furthermore, under specific assumptions, the generating function is the eigenvalue distribution symbol as well \cite{Barbarino2020-eo,bottcher051,Golinskii2007-ue,Tilli1999-fk}. Further spectral results of extremal type are in \cite{bottcher981,capizzano991}, while we refer to the book by Garoni and Serra-Capizzano \cite{garoni171} for a general account on Toeplitz-like matrix-sequences, variable coefficient generalizations, and several applications.

In this article we focus on the spectrum of Toeplitz matrices generated by symbols of the form 
\begin{align}
f(\theta)=\hat{f}_0+\hat{f}_{r}e^{\mathbf{i}r\theta}+\hat{f}_{-s}e^{-\mathbf{i}s\theta},\label{eq:symbol}
\end{align}
where \(\theta\in\left[-\pi, \pi\right]\), \(r,s\in\mathbb{Z^+}\), \(r\neq s\), \(\hat{f}_r, \hat{f}_0, \hat{f}_{-s} \in \mathbb{C}\) and \(\hat{f}_{r}, \hat{f}_{-s} \neq 0\). Without loss of generality we take $r < s$ (since, if $r > s$, we can instead study the transposed matrix generated by \(f(-\theta)\)).
The generated Toeplitz matrices $T_n(f)$ are of the form
\begin{align*}
T_n(f)=
\left[
\begin{array}{cccccccccc}
\hat{f}_0&&&\hat{f}_{-s}&&\\
&\hat{f}_0&&&\hat{f}_{-s}&\\
\hat{f}_{r}&&\ddots&&&\ddots\\
&\hat{f}_{r}&&\ddots&&&\ddots\\
&&\ddots&&\ddots&&&\hat{f}_{-s}\\
&&&\ddots&&\ddots&\\
&&&&\ddots&&\ddots&\\
&&&&&\hat{f}_{r}&&\hat{f}_0
\end{array}
\right].
\end{align*}
We can define the ``symmetrised'' symbol
\begin{align}
    \Bar{f}(\theta)&=\hat{f}_0 + \hat{f}_{r}^{s/(r+s)}\hat{f}_{-s}^{r/(r+s)} \left(e^{\mathbf{i}r\theta}+e^{-\mathbf{i}s\theta}\right),\label{eq:symmetrised}
\end{align}
and by a similarity transformation we can see that \(T_n(f)\sim T_n(\Bar{f})\).
We note the following two cases for which results are already known,
\begin{itemize}
    \item \(r=s=1\): the classical tridiagonal Toeplitz matrices; see, e.g.,~\cite{bottcher051,bottcher001, krein501, trefethen051,V_Egervary1928-cx}, and the references therein,
    \item \(r=s>1\): the ``symmetrically sparse tridiagonal'' Toeplitz matrices; see, e.g., ~\cite{ekstrom181,Ekstrom2020-me,Fonseca2020-yx}, and the references therein.
\end{itemize}
\begin{rmrk}
For completeness, we also mention the trivial case where \(r<0\) or \(s<0\), that is, the two off diagonals are on the same side of the main diagonal. Then, all the eigenvalues of \(T_n(f)\) are equal to \(\hat{f}_0\).
\end{rmrk}
\noindent
This paper focuses on the eigenvalues of the more general case of \(r\neq s\). To that end, we can simplify the problem by noting that for any \(\{r,s\}\) the eigenvalues \(\lambda_j(T_n(f))\), for \(j=1,\ldots,n\), can be written as
\begin{align}
\lambda_j(T_n(f))=\lambda_j(T_n(\Bar{f}))=\hat{f}_0+\hat{f}_{r}^{s/(r+s)}\hat{f}_{-s}^{r/(r+s)}\lambda_j(T_n(g_{r,s})),\qquad j = 1,\dots, n,\label{intro:spectrum_reduced}
\end{align}
where \(g_{r,s}\) is the generating symbol
\begin{align}
    g_{r,s}(\theta) = e^{\mathbf{i}r\theta}+e^{-\mathbf{i}s\theta}.
    \label{eq:intro:g}
\end{align}
The symbol \(g_{r,s}(\theta)\) generates the matrices
\begin{align}
    T_n(g_{r,s})=
    \left[
    \begin{array}{cccccccccc}
    &&&1&&\\
    &&&&1&\\
    1&&&&&\ddots\\
    &1&&&&&\ddots\\
    &&\ddots&&&&&1\\
    &&&\ddots&&&\\
    &&&&\ddots&&&\\
    &&&&&1&&
    \end{array}
    \right],
    \label{eq:intro:Tng}
\end{align}
where we have ones on the \(r\)-th and \(-s\)-th diagonals. We may without loss of generality restrict our study to the generating symbol \(g_{r,s}(\theta)\), since the spectrum of \(T_n(f)\) may be recovered from the spectrum of \(T_n(g_{r,s})\) through shifting, rotation, and scaling by \eqref{intro:spectrum_reduced}.

The article is organised as follows. In Section~\ref{sec:main} we first propose Conjecture~\ref{conj:1} regarding the eigenvalues of \(T_n(g_{r,s})\) in \eqref{eq:intro:Tng} (where \(g_{r,s}(\theta)\) is defined in \eqref{eq:intro:g}).  
Then, the properties of the full spectrum of \(T_n(g_{r,s})\) is briefly discussed in subsection~\ref{sec:main:propfull}, and  the positive real eigenvalues in subsection~\ref{sec:main:propposreal}. In subsection~\ref{sec:main:constructposreal} it is proposed how the positive real spectrum of \(T_n(g_{r,s})\) may be computed. Subsequently, in subsection~\ref{sec:main:fullspectrum} it is shown how to construct the full spectrum of \(T_n(g_{r,s})\), from the positive real spectrum. In Section~\ref{sec:examples} examples and numerical experiments, supporting Conjecture~\ref{conj:1} are presented. In Section~\ref{sec:conclusions} open problems, conclusions,  and future avenues of research are presented. Finally, in Appendix~\ref{sec:appendix} full code in \textsc{Julia}~\cite{Julia-2017} is supplied for the proposed algorithm.

\section{Main Results}
\label{sec:main}

The main result of this article is the following conjecture, with the rest of the article dedicated to supporting it.
\begin{conj}[The eigenvalues of Toeplitz matrices generated by \(g_{r,s}(\theta)\)]
    \label{conj:1}
Assume \(g_{r,s}(\theta)\) is defined as in \eqref{eq:intro:g} (with \(1\leq r\leq s\)). Then, we know that a subset of the eigenvalues \(\lambda_j(T_n(g_{r,s}))\), for \(j=1,\ldots, n\), are real and positive~\cite{schmidt601}; denote this subset of the eigenvalues as \(\lambda_+(T_n(g_{r,s}))\).
Then, the eigenvalues in \(\lambda_+(T_n(g_{r,s}))\) are, for \(\gamma=\mathrm{gcd}(r,s)\), given by 
\begin{align}
    \lambda_+ \left(T_n(g_{r,s})\right)^{\omega} &= \left(\bigcup_{k = 0}^{\gamma - \beta_\gamma}\lambda_+\left(T_{n_\gamma}(g_{r_\gamma,s_\gamma})\right)\right) \bigcup \left(\bigcup_{k = 0}^{\beta_\gamma}\lambda_+\left(T_{n_\gamma+1}(g_{r_\gamma,s_\gamma})\right)\right)
\end{align} 
where the parameters \(\gamma, \omega, \beta_\gamma, n_\gamma\), and \(r_\gamma\), \(s_\gamma\), are defined in \eqref{eq:gamma}, \eqref{eq:omegas}, \eqref{eq:bg}, \eqref{eq:ng}, and \eqref{eq:rgsg}. The eigenvalues \(\lambda_+\left(T_{n_\gamma}(g_{r_\gamma,s_\gamma})\right)\) (and \(\lambda_+\left(T_{n_\gamma+1}(g_{r_\gamma,s_\gamma})\right)\)) are exactly the same as those of a smaller matrix \(B_{(n_\gamma)_\sigma}^{n_\gamma,r_\gamma,s_\gamma}\) (and \(B_{(n_\gamma+1)_\sigma}^{n_\gamma+1,r_\gamma,s_\gamma}\)), possibly with multiplicity; a construction of these matrices is presented in Section~\ref{sec:main:constructposreal}.

The full spectrum of \(T_n(g_{r,s})\) can be constructed by \(\omega\) rotations of \(\lambda_+(T_n(g_{r,s}))\) plus \(n_0\) zero eigenvalues; see Section~\ref{sec:main:fullspectrum} and \eqref{eq:n0}.
\end{conj}

Hence, to summarise, we claim in Conjecture~\ref{conj:1} that the positive real eigenvalues of a matrix \(T_n(g_{r,s})\), where \(1\leq r\leq s\), are the same as the \(\omega\):th root of the eigenvalues of a non-unique other matrix (or two, with specified multiplicity) which we can construct automatically. The full spectrum of \(T_n(g_{r,s})\) can then be constructed by rotations of these positive real eigenvalues, plus possibly \(n_0\) zeros.

\subsection{Properties of \(\lambda_j(T_n(g_{r,s}))\)  for \(j=1,\ldots,n\)}
\label{sec:main:propfull}
\begin{figure}[H]
\centering
\includegraphics[width=0.32\textwidth]{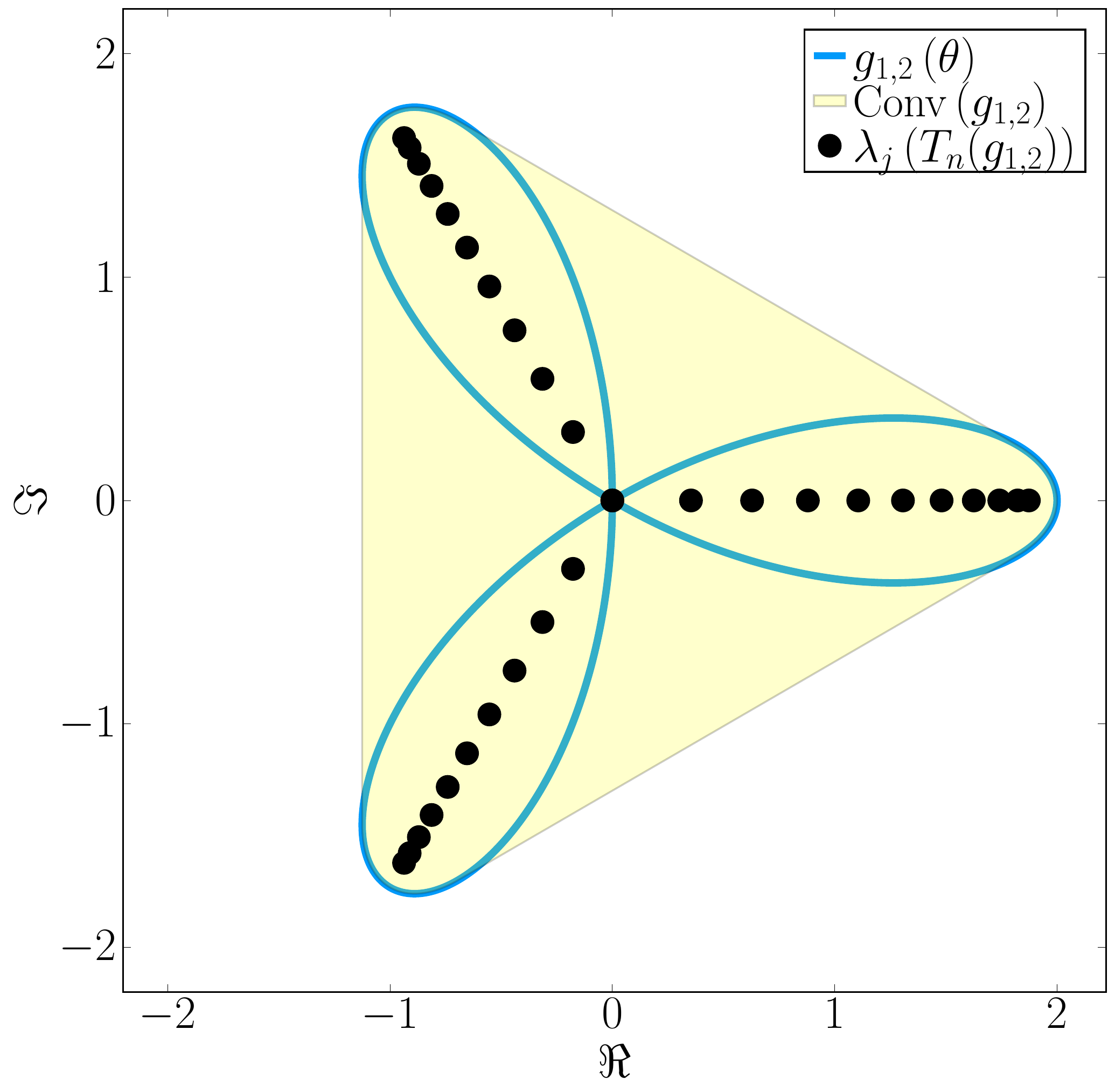}
\includegraphics[width=0.32\textwidth]{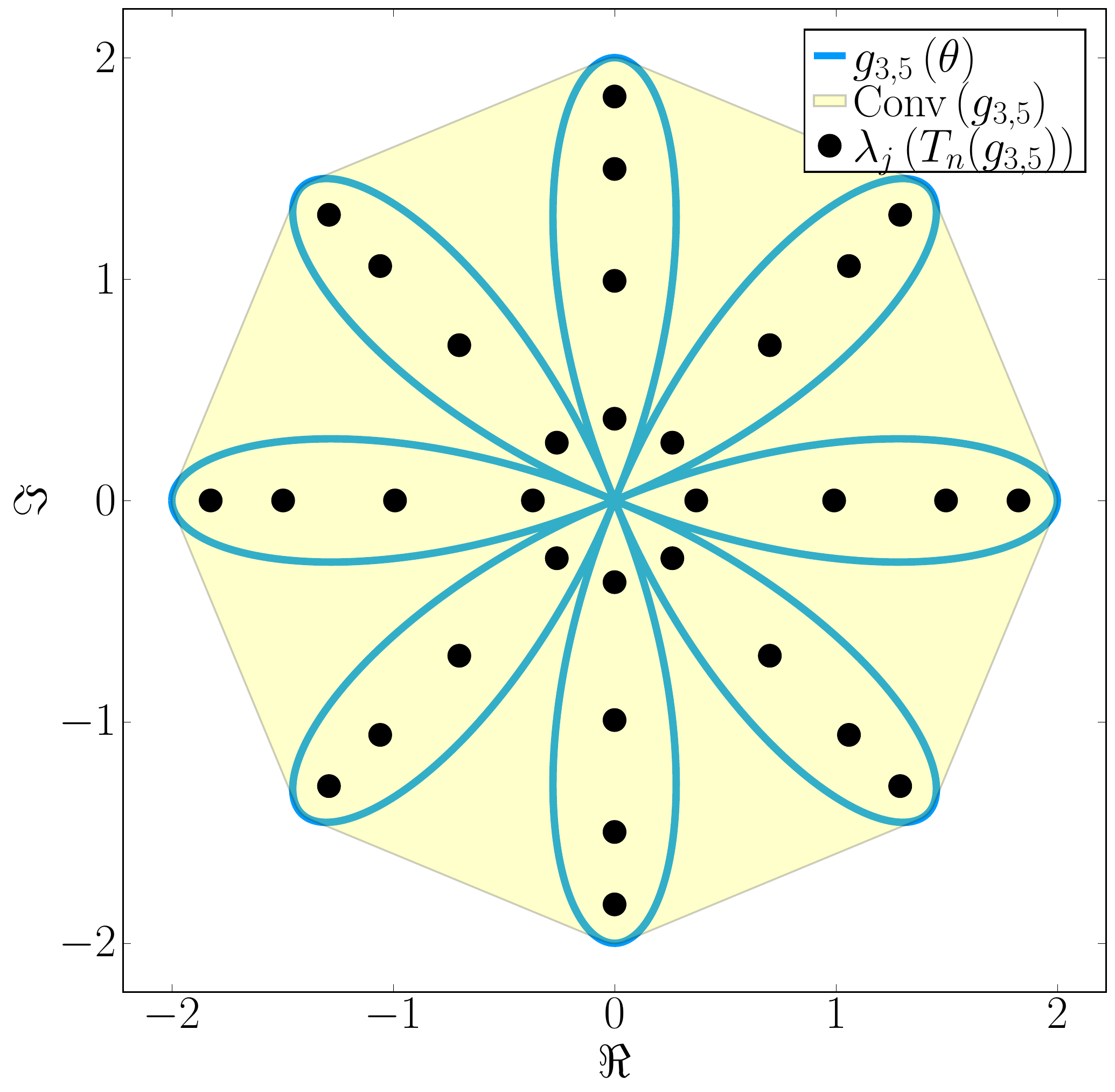}
\includegraphics[width=0.32\textwidth]{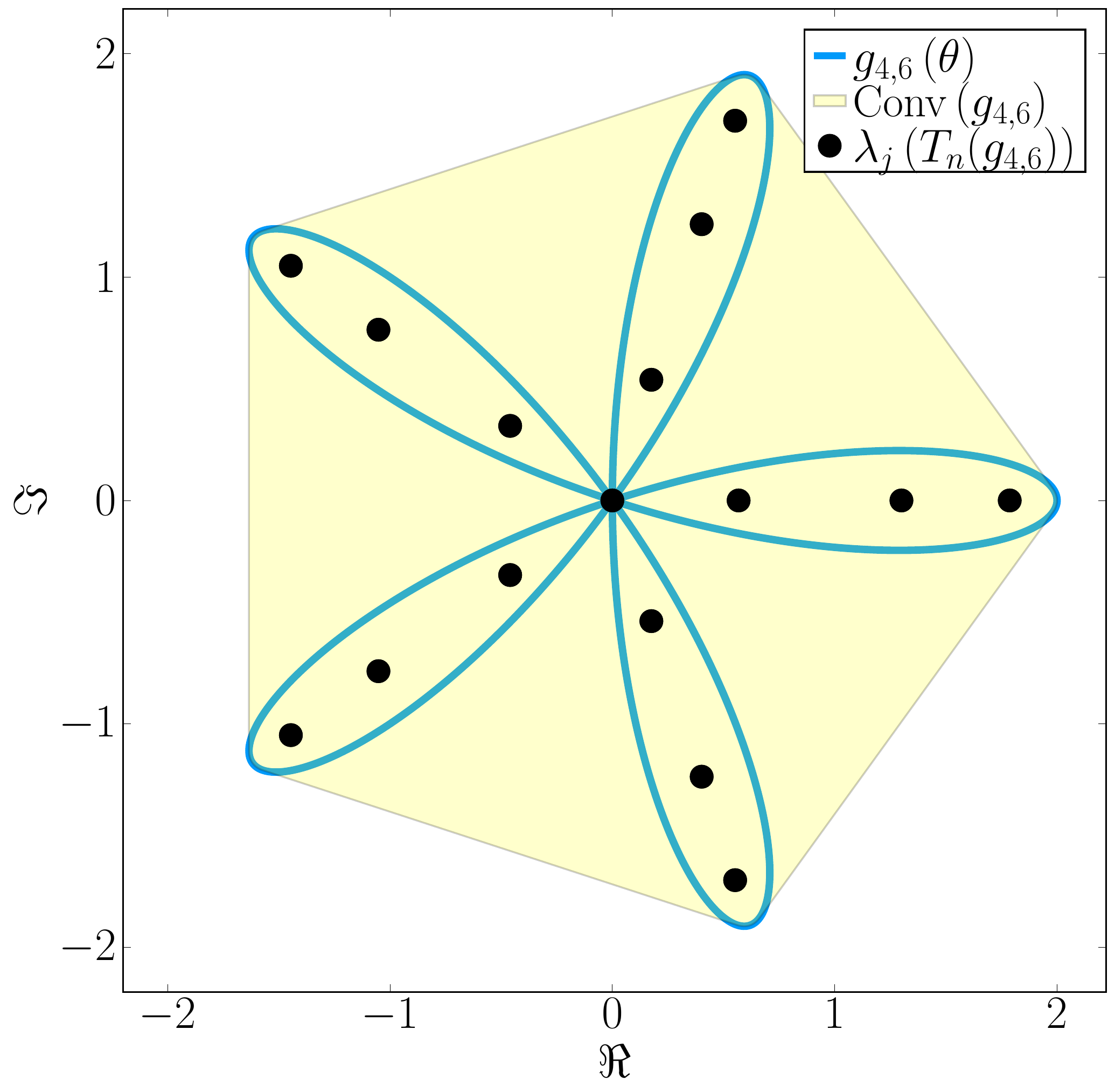}
\caption{Symbol \(g_{r,s}\left(\theta\right)\), the convex hull Conv\(\left(g_{r,s}\right)\) and the eigenvalues \(\lambda_j(T_n(g_{r,s}))\), for \(j = 1,\dots,n\). Parameters \(\{n,r,s\}\). \textbf{Left:} \(\{32, 1, 2\}\). \textbf{Middle:} \(\{32, 3, 5\}\). \textbf{Right:} \(\{32, 4, 6\}\).}
\label{fig:symbolsspectra}
\end{figure}

\begin{figure}[H]
    \centering
    \includegraphics[width=0.32\textwidth]{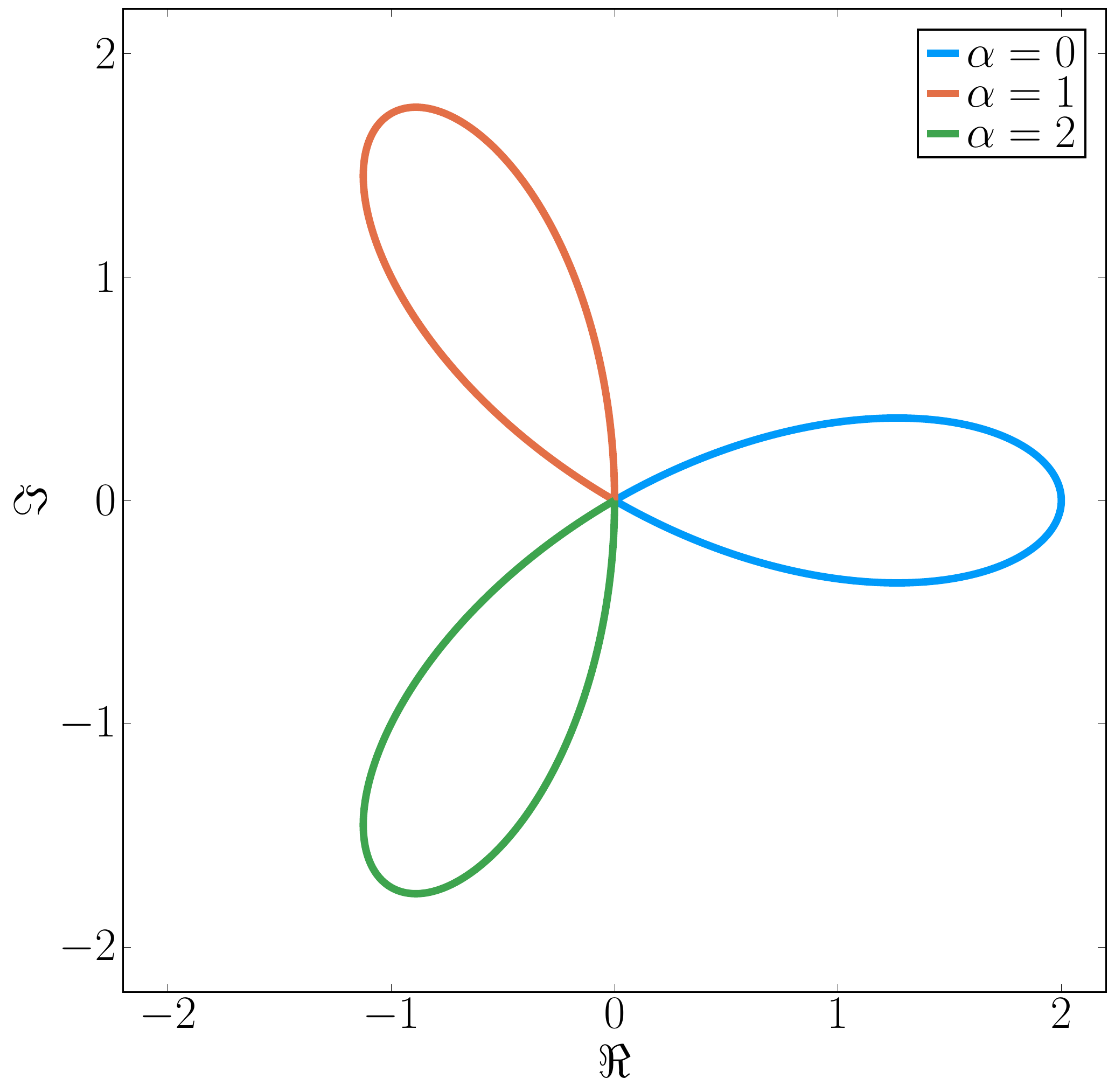}
    \includegraphics[width=0.32\textwidth]{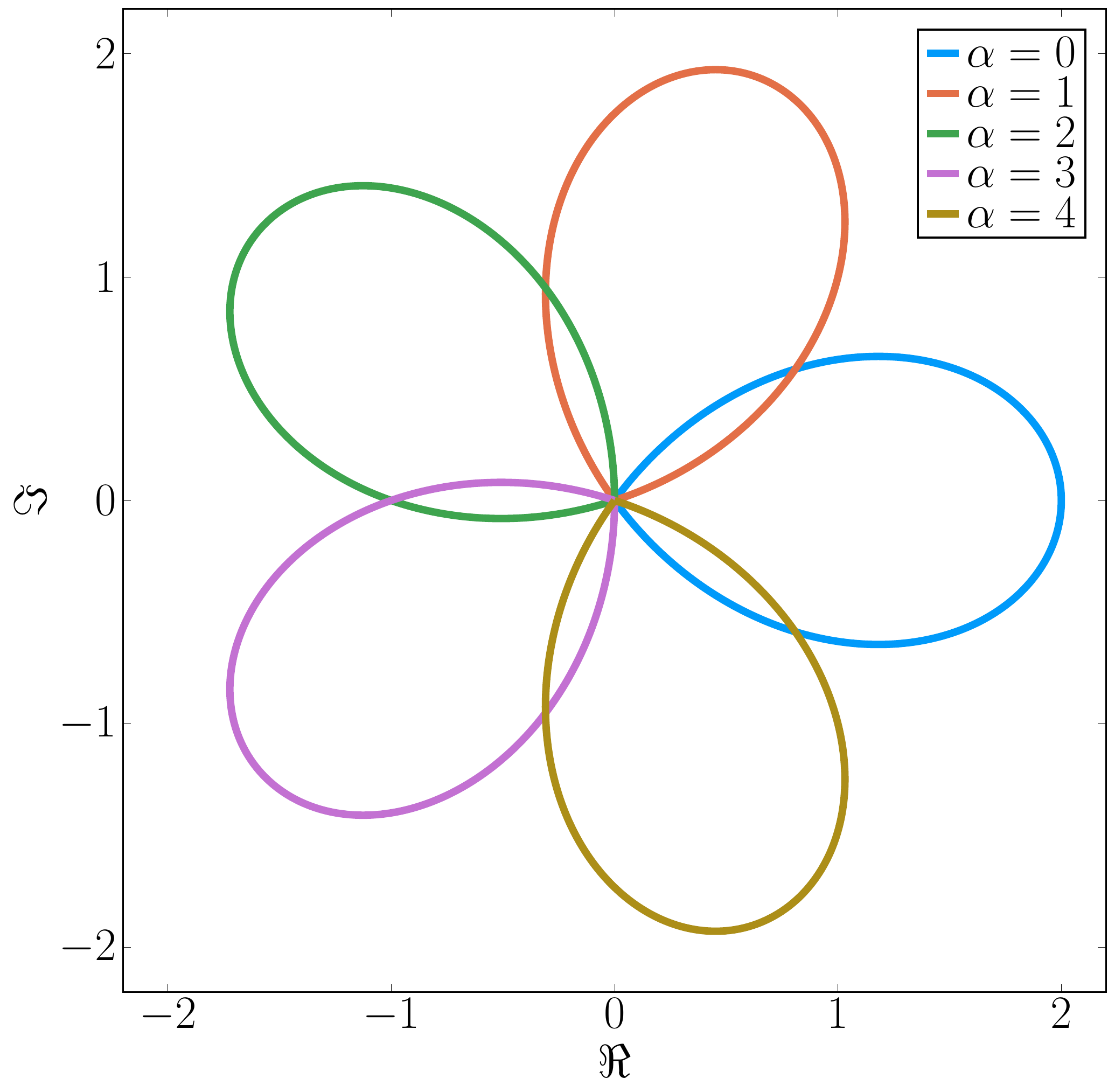}
    \includegraphics[width=0.32\textwidth]{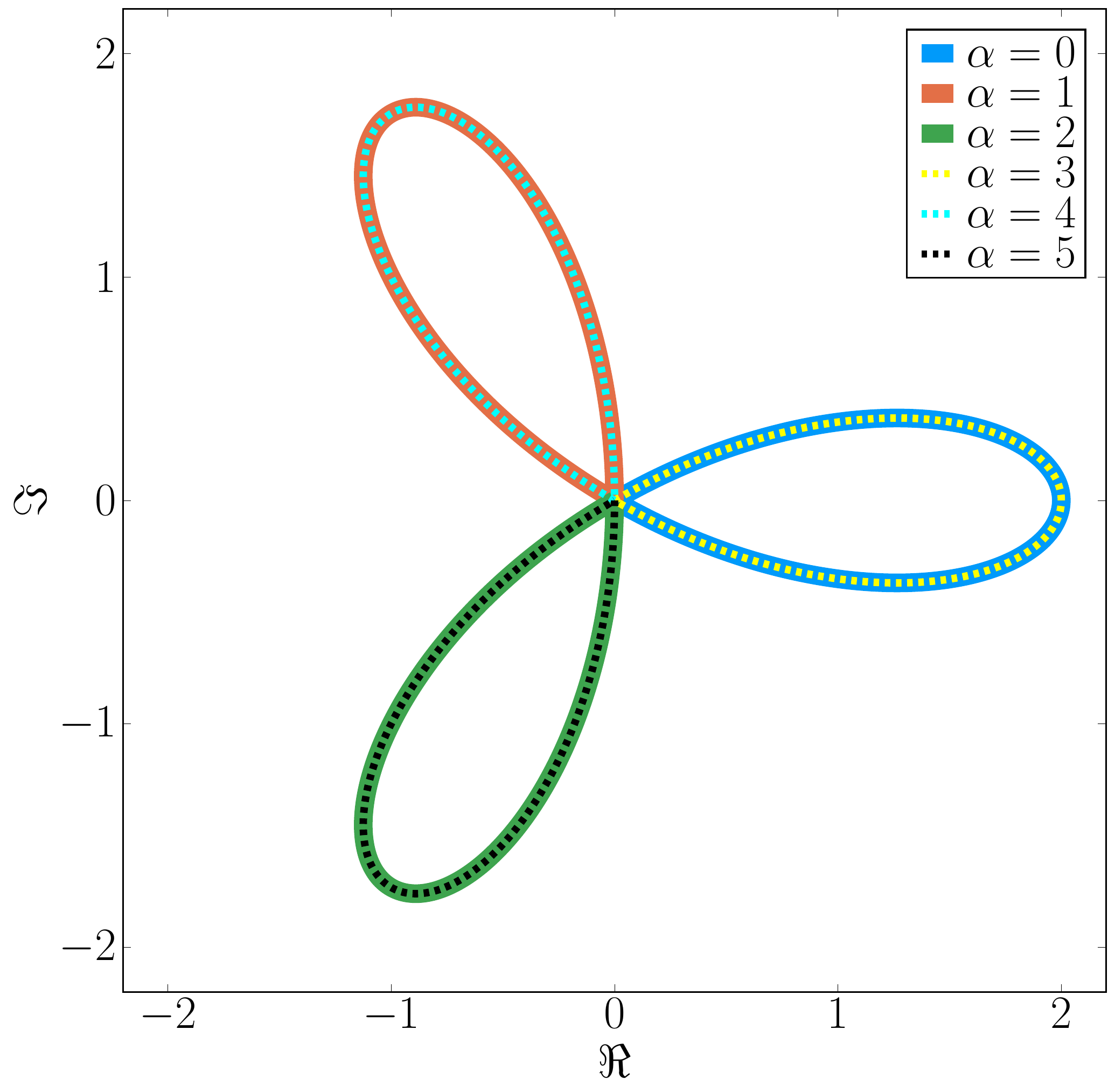}
    \caption{Symbol \(g_{r,s}\left(\theta\right)\), with the different arms \(\alpha=0,\ldots,\sigma-1\) highlighted (number of distinct arms is \(\omega\)). Parameters \(\{r,s\}\). 
    \textbf{Left:} \(\{1, 2\}\) where \(\omega=\sigma=3\). 
    \textbf{Middle:} \(\{1, 4\}\) where \(\omega=\sigma=5\). 
    \textbf{Right:} \(\{2, 4\}\) where \(\omega=3\) and \(\sigma=6\).}
    \label{fig:symbolsarms}
\end{figure}

It is known that the eigenvalues of \(T_n(g_{r,s})\), denoted \(\lambda_j\left(T_n(g_{r,s})\right)\) for \(j=1,\ldots, n\), lie in the convex hull of the essential range of the generating symbol \(g_{r,s}(\theta)\); see, e.g.,~\cite{bottcher981,bottcher001,capizzano991,garoni171}, and the references therein. 
This fact is visualised in Figure~\ref{fig:symbolsspectra}, for three different parameters \(\{n,r,s\}\).
The symbols \(g_{r,s}\), and thus the eigenvalues \(\lambda_j(T_n(g_{r,s}))\) for \(j=1,\ldots,n\) lie in a ``star'' formation, for example see~\cite[pp. 204--205, p. 275]{bottcher051}, with \(\omega\) distinct ``arms'', where
\begin{align}
\gamma&=\mathrm{gcd}(r,s),\label{eq:gamma}\\
\sigma &= r + s,\nonumber\\
\omega&=\frac{\sigma}{\gamma},\label{eq:omegas}\\
\beta_\sigma&=\mathrm{mod}(n,\sigma),\label{eq:betas}\\
n_\sigma&=\frac{n-\beta_\sigma}{\sigma}.\label{eq:ns}
\end{align}
The matrix \(T_n(g_{r,s})\) has a total of
\begin{align}
    n_0&=(\gamma - \beta_\gamma)\mathrm{mod}(n_\gamma,\omega) + \beta_\gamma\mathrm{mod}(n_\gamma+1,\omega),\label{eq:n0}
\end{align}
zero eigenvalues, and in each distinct arm we have a total of \(\frac{n-n_0}{\omega}\) non-zero eigenvalues. In the three panels of Figure~\ref{fig:symbolsspectra} we have respectively two, none, and two zero eigenvalues.
Generally we have that arm \(\alpha\) is generated by \(g_{r,s}\left(\frac{\theta}{\sigma} + \alpha\frac{2\pi}{\sigma}\right)\) where \(\alpha=0\) corresponds to the arm that covers the positive real axis~\cite{schmidt601}, \(\alpha=1\) is the next arm anti-clockwise, etc.; see Figure~\ref{fig:symbolsarms} for a visual reference.

\subsection{Properties of \(\lambda_+(T_n(g_{r,s}))\)}
\label{sec:main:propposreal}
It is reported in~\cite{bottcher051,schmidt601}, referring to results in~\cite{biernacki281}, that there is one arm of the star, formed by the symbol \(g_{r,s}\), that covers a part of the positive real axis.  The eigenvalues \(\lambda_j(T_n(g_{r,s}))\) in this arm are positive and real-valued.
Denote the set of positive real-valued eigenvalues \(\lambda_+(T_n(g_{r,s}))\), and we know from~\cite{schmidt601} that
\begin{align*}
0<\lambda_+(T_n(g_{r,s}))\leq R=\frac{\sigma}{r^{\frac{r}{\sigma}}s^{\frac{s}{\sigma}}}.
\end{align*}
Note that we here exclude the zero eigenvalues, which are of multiplicity \(n_0\), defined in \eqref{eq:n0}. 
The non-zero eigenvalues of the other arms are rotations of \(\lambda_+(T_n(g_{r,s}))\); see~\cite{schmidt601}.
Now define the following function for the arm covering the eigenvalues \(\lambda_+(T_n(g_{r,s}))\)
\begin{align*}
    a_{r,s}(\theta) = g_{r,s}\left(\frac{\theta}{\sigma}\right) = e^{\mathbf{i}\frac{r}{\sigma}\theta}+e^{-\mathbf{i}\frac{s}{\sigma}\theta}.
\end{align*}
We notice that we do not have integer exponents, but if we construct the following symbol from \(a_{r,s}(\theta)\),
\begin{align}
    b_{r,s}(\theta) = a_{r,s}(\theta)^\sigma = e^{\mathbf{i}r\theta}\left(1 + e^{-\mathbf{i}\theta}\right)^\sigma,\label{eq:b}
\end{align}
we have a standard symbol for generated Toeplitz or Toeplitz-like matrices. 
\begin{rmrk}
We note that Toeplitz matrices \(T_n(b_{r,s})\), generated by the symbol \(b_{r,s}(\theta)\) in \eqref{eq:b}, are non-symmetric, but have only real eigenvalues~\cite{Shapiro2019}. That is, the symbol \(b_{r,s}\) is complex valued, and we have \(\{T_n(b_{r,s})\}\sim_{\textsc{glt},\sigma}b_{r,s}\) and \(\{T_n(b_{r,s})\}\sim_{\lambda}\mathfrak{b}_{r,s}\), where \(\mathfrak{b}_{r,s}\) is real valued, see \cite[Section 4.3]{Shapiro2019},
\begin{align}
\mathfrak{b}_{r,s}(\theta)=\frac{\sin^\sigma(\theta)}{\sin^r(\frac{r}{\sigma}\theta)\sin^s(\frac{s}{\sigma}\theta)}.\label{eq:shapirob}
\end{align}
See the theory of generalised locally Toeplitz (GLT) sequences~\cite{garoni171} where we can deduce that also Toeplitz-like matrices \(B_n\) with the symbol \(b_{r,s}(\theta)\) are \(\{B_n\}\sim_{\textsc{glt},\sigma}b_{r,s}\) and \(\{B_n\}\sim_{\lambda}\mathfrak{b}_{r,s}\) (except for possibly \(o(n)\) outliers, which are not present in the current setting of this article).
\end{rmrk}
\noindent Furthermore, it is possible to rewrite the symbol \(b_{r,s}(\theta)\) in \eqref{eq:b} as a product of symbols
\begin{align}
    b_{r,s}(\theta) = c(\theta)^{r}c(-\theta)^{s}
    \label{eq:bcmc}
\end{align}
where
\begin{align*}
    c(\theta) = 1 + e^{\mathbf{i}\theta}.
\end{align*}
The fundamental idea of Conjecture~\ref{conj:1} is that we can construct a matrix with symbol \(b_{r,s}\) (or two with symbol \(b_{r_\gamma,s_\gamma}\)) of size smaller than \(n\), whose eigenvalues (possibly with multiplicity) are exactly the eigenvalues \(\lambda_+(T_n(g_{r,s}))\) to the power of \(\omega\).

\subsection{Computing \(\lambda_+(T_n(g_{r,s}))\)}
\label{sec:main:constructposreal}

We now detail how to compute the positive real eigenvalues of \(T_n(g_{r,s})\). In the case \(\gamma = 1\), many parts of the ideas described below will simplify, however we here present the general case of \(\gamma > 1\). When \(\gamma > 1\), we will have a number of non-zero eigenvalues with multiplicity larger than one, and we can then look at a ``reduced'' case which yields the eigenvalues we seek. First define the following variables
\begin{align}
    \beta_\gamma &= \mathrm{mod}(n, \gamma),\label{eq:bg}\\
    n_\gamma &= \frac{n - \beta_\gamma}{\gamma},\label{eq:ng}\\
    r_\gamma &= \frac{r}{\gamma}, \quad s_\gamma = \frac{s}{\gamma}\label{eq:rgsg}\\
    \sigma_\gamma&=r_\gamma+s_\gamma=\omega.\label{eq:sigmag}
\end{align}
Then, we use \(\{n_\gamma, r_\gamma, s_\gamma\}\) as the new triple, that is, find the eigenvalues \(\lambda_+(T_{n_\gamma}(g_{r_\gamma,s_\gamma}))\), by applying the algorithm as described in Section~\ref{sec:main:bg1}. 
We use the following notation for variables generated for the construction of the matrix \(B_{(n_\gamma)_\sigma}^{n_\gamma,r_\gamma,s_\gamma}\)
\begin{align}
    (\beta_\gamma)_\sigma&=\mathrm{mod}(n_\gamma,\sigma_\gamma),\label{eq:betags}\\
    (n_\gamma)_\sigma&=\frac{n_\gamma-(\beta_\gamma)_\sigma}{\sigma_\gamma}.\label{eq:ngs}
\end{align}
After computing the positive real eigenvalues associated with this triple, we then repeat all the eigenvalues \(\gamma - \beta_\gamma\) times.

However, if \(\beta_\gamma > 0\) we also need to compute the eigenvalues for the triple \(\{n_\gamma+1, r_\gamma, s_\gamma\}\), that is, find the eigenvalues \(\lambda_+(T_{n_\gamma+1}(g_{r_\gamma,s_\gamma}))\) where we compute the values in (\ref{eq:betags}) and (\ref{eq:ngs}) using \(n_\gamma + 1\) instead of \(n_\gamma\). These eigenvalues are repeated \(\beta_\gamma\) times. 

\noindent
Thus, the full computation of \(\lambda_+ \left(T_n(g_{r,s})\right)\) is given by
\begin{align*}
    \lambda_+ \left(T_n(g_{r,s})\right)^{\omega} &= \left(\bigcup_{k = 0}^{\gamma - \beta_\gamma}\lambda_+\left(T_{n_\gamma}(g_{r_\gamma,s_\gamma})\right)\right) \bigcup \left(\bigcup_{k = 0}^{\beta_\gamma}\lambda_+\left(T_{n_\gamma+1}(g_{r_\gamma,s_\gamma})\right)\right).
\end{align*}
We will now presnt how to construct the matrices which has the same eigenvalues as \(\lambda_+(T_{n_\gamma}(g_{r_\gamma,s_\gamma}))\), \(\lambda_+(T_{n_\gamma + 1}(g_{r_\gamma,s_\gamma}))\) or even all of \(\lambda_+(T_{n}(g_{r,s}))\) in the case \(\gamma = 1\). Thus, we may now assume that going forward the \(\gcd(r,s) = 1\), whether those be the original \(\{r, s\}\) if \(\gamma = 1\) or \(\{r_\gamma,s_\gamma\}\) if we are working with the ``reduced'' case where \(\gamma > 1\).

\subsubsection{Construction of matrix \(B_{n_\sigma}^{n,r,s}\)}
\label{sec:main:bg1}

We propose that there exists a Toeplitz-like matrix, which we denote \(B_{n_\sigma}^{n,r,s}\), which has the symbol \(b_{r,s}(\theta)\) and has the eigenvalues \(\cup_{j=1}^{n_\sigma}(\lambda_j(B_{n_\sigma}^{n,r,s}))^\frac{1}{\sigma} = \lambda_+(T_{n}(g_{r,s}))\).
The matrix \(B_{n_\sigma}^{n,r,s}\) is here uniquely defined by \(\{n,r,s\}\) in the algorithm below (however, it is non-unique in sharing the eigenvalues with \(\lambda_+(T_n(g_{r,s}))\), and can be constructed in many ways). The subscript \(n_\sigma\), defined in \eqref{eq:ns}, is only explicitly written to emphasise the size of the matrix; \(n_\sigma<n\).

To construct a matrix \(B_{n_\sigma}^{n,r,s}\) we use a prescribed, but non-unique, product of matrices of the form \(T_{n_\sigma}(e^{-\mathbf{i}\theta}c(\theta)^m)\) and their transposes; an implementation to generate these matrices is presented in Appendix~\ref{code:tnc}.
\begin{rmrk}
The non-zero Fourier coefficients \(\hat{f}_k\), \(k=-1,\ldots,m-1\) of the symbol \(e^{-\mathbf{i}\theta}c(\theta)^m=e^{-\mathbf{i}\theta}(1+e^{\mathbf{i}\theta})^m\), used to generate the matrices \(T_{n_\sigma}(e^{-\mathbf{i}\theta}c(\theta)^m)\), are given by
\begin{align*}
\hat{f}_k=\binom{m}{k+1}.
\end{align*}
\end{rmrk}
We provide here an algorithm to automatically generate the matrix \(B_{n_\sigma}^{n,r,s}\). We begin by constructing two matrices \(\mathcal{M}, \mathcal{P}\in \mathbb{Z}^{\sigma\times r}\) (defined in (\ref{results:Mmatrix}) and (\ref{results:Pmatrix})). The elements of these matrices define how to generate \(B_{n_\sigma}^{n,r,s}\), 
\begin{align}
    B_{n_\sigma}^{n,r,s} 
    &= \prod_{k=1}^{r} T_{n_\sigma}^{\mathop\intercal}\left(e^{-\mathbf{i}\theta}c(\theta)^{m_{k}}\right) \left(T_{n_\sigma}(e^{-\mathbf{i}\theta}c(\theta))\right)^{p_{k}},\label{eq:realEigMatCode}\\
    &= \prod_{k=1}^{r} T_{n_\sigma}\left(c(\theta)c(-\theta)^{m_{k}-1}\right) \left(T_{n_\sigma}(c(-\theta))\right)^{p_{k}},\label{eq:realEigMat}
\end{align}
where \(m_k = \mathcal{M}_{\beta+1, k}\) and \(p_k = \mathcal{P}_{\beta+1, k}\). See Appendix~\ref{code:bns} for an implementation of \eqref{eq:realEigMatCode}. Furthermore, we have that the symbol of \( B_{n_\sigma}^{n,r,s}\) in \eqref{eq:realEigMat} will be 
\begin{align*}
\{ B_{n_\sigma}^{n,r,s}\}_{n_\sigma}&\sim_\sigma \prod_{k=1}^r c(\theta)c(-\theta)^{m_{k}-1}c(-\theta)^{p_k}\\
&=c(\theta)^rc(-\theta)^{\sum_{k=1}^r(m_k+p_k)-r}\\
&=\left[\sum_{k=1}^r(m_k+p_k)=\sigma \text{ by construction}\right]\\
&=c(\theta)^rc(-\theta)^s\\
&=b_{r,s}(\theta)
\end{align*}
as expected from~\eqref{eq:bcmc}. Note that all \(\{ B_{n_\sigma}^{n,r,s}\}_{n_\sigma}\sim_\lambda \mathfrak{b}_{r,s}\).

\subsubsection*{Construction of matrices \(\mathcal{M}\) and \(\mathcal{P}\)}
\label{sec:main:mp}
Defining \(\tau=\mathrm{mod}(s,r)\), then the elements of \(\mathcal{M}\) and \(\mathcal{P}\) are given by the following construction.
The matrices \(\mathcal{M}, \mathcal{P}\in \mathbb{Z}^{\sigma\times r}\) are both rectangular Toeplitz matrices defined by
\begin{align}
    \mathcal{M}_{i,j} &= \begin{cases}
        1, & \text{if } i = 1,\\
        \mathcal{M}_{i-1,j}+1, & \text{if } i > 1 \text{ and } \mathrm{mod}(j-i,r) = \mathrm{mod}(-1,r),\\
        \mathcal{M}_{i-1,j}, & \text{otherwise},
        \end{cases}\label{results:Mmatrix}\\
        \mathcal{P}_{i,j} &= \begin{cases}
            \frac{s-\tau}{r}, & \text{if } i = 1 \text{ and } j \leq r -  \tau,\\
            \frac{s-\tau}{r} + 1, & \text{if } i = 1 \text{ and } j > r -  \tau, \\
            \mathcal{P}_{i-1,j}-1, & \text{if } i > 1 \text{ and } \mathrm{mod}(j-i, r) = \mathrm{mod}(r-\tau-1, r),\\
            \mathcal{P}_{i-1,j}, & \text{otherwise}.
            \end{cases}\label{results:Pmatrix}
\end{align}
An alternative, but equivalent, construction of \(\mathcal{M}\) and \(\mathcal{P}\) is presented in Appendix~\ref{code:MP}.
In the case that \(r=1\) then the matrices \(\mathcal{M}\) and \(\mathcal{P}\) are given as above, however if \(r>1\), then the columns of the matrices \(\mathcal{M}\) and \(\mathcal{P}\) in \eqref{results:Mmatrix} and \eqref{results:Pmatrix} need to be permuted and so will no longer have a Toeplitz structure. We create the permutation by first creating the ordered set \(\mathcal{S}_\tau = \left(\tau, 2\tau,3\tau, \ldots, r\tau\right)\). Then, the permutation is defined as \(p_{\mathrm{perm}} = \mathrm{mod}(\mathcal{S}_\tau, r)\) where the modulo operator is applied element-wise to the ordered set \(\mathcal{S}_\tau\), and instead of 0 for the final element (since it will be \(\bmod(r\tau, r) = 0\)) we will use \(r\). Then, permute the columns of \(\mathcal{M}\) and \(\mathcal{P}\) with \(p_{\mathrm{perm}}\). See an implementation in Appendix~\ref{code:MP}.

\subsubsection*{Erroneous perturbations when constructing \(B_{n_\sigma}^{n,r,s}\)}

For \(\beta_\sigma > s\) we will have that for some \(k\) in \eqref{eq:realEigMat},~\(p_{k}\) are negative; that is, inverses of \(T_{n_\sigma}(c(-\theta))\) are used in the construction of \( B_{n_\sigma}^{n,r,s}\). This results in incorrect and unwanted non-zero integer-valued perturbations being generated in the top right corner due to the inverted matrices; the generated \(B_{n_\sigma}^{n,r,s}\) is in fact \( B_{n_\sigma}^{n,r,s}+R_{n_\sigma}^{n,r,s}\) for odd \(n_\sigma\) and  \( B_{n_\sigma}^{n,r,s}-R_{n_\sigma}^{n,r,s}\) for even \(n_\sigma\), where \(R_{n_\sigma}^{n,r,s}\) is a low rank perturbation with strictly positive entries.

Our current approach to address this is to subtract these perturbations from the matrix and then proceed normally; that is, identify \(R_{n_\sigma}^{n,r,s}\) and remove it. If the perturbations lie within the non-zero bandwidth of the correct matrix, we can construct a larger matrix where the perturbations are outside the bandwidth of the matrix, identify \(R_{n_\sigma}^{n,r,s}\), and add or subtract it from the corresponding elements in the original sized (\(n_\sigma\)) matrix. See an implementation in Appendix~\ref{code:bns}. However, this does not work in the case that the unwanted values cross over the main diagonal, and so we then have the restriction that \(n>(r-1)\sigma\) rather than the essential restriction \(n> s\). It is worth noting that by an informed search, a correct matrix  \( B_{n_\sigma}^{n,r,s}\) for \(s< n \leq (r-1)\sigma\) can always manually be found, but an efficient algorithm to generate them has not yet been designed. We reiterate that this restriction on \(n\) is only necessary in the case that there are erroneous corner elements.

To exemplify these erroneous elements given by \(R_{n_\sigma}^{n,r,s}\), that have to be removed when generating \( B_{n_\sigma}^{n,r,s}\), see Section~\ref{sec:example:CornerAndPerms}.

\subsection{Constructing \(\lambda_j(T_n(g_{r,s}))\) for \(j=1,\ldots,n\) from \(\lambda_+(T_n(g_{r,s}))\)}
\label{sec:main:fullspectrum}
Finally, we can then construct all the eigenvalues of \(T_n(g_{r,s})\) from \(\lambda_+(T_n(g_{r,s}))\) by
\begin{align*}
    \lambda\left(T_n(g_{r,s})\right) &= \left(\bigcup_{\alpha = 0}^{\omega-1} e^{2\pi\mathbf{i}\frac{\alpha}{\omega}} \lambda_+ \left(T_n(g_{r,s})\right)\right) \cup \left(\mathbf{0}_{n_0}\right)\\
    &=\left(\left[e^{2\pi\mathbf{i}\frac{0}{\omega}}, e^{2\pi\mathbf{i}\frac{1}{\omega}},\ldots, e^{2\pi\mathbf{i}\frac{\omega-1}{\omega}}\right]\otimes \lambda_+ \left(T_n(g_{r,s})\right)\right) \cup \left(\mathbf{0}_{n_0}\right)
\end{align*}
where \(\omega\) is as defined in (\refeq{eq:omegas}), the multiplication of \(e^{2\pi\mathbf{i}\frac{\alpha}{\omega}}\), for \(\alpha=0,\ldots,\omega-1\), with \(\lambda_+ \left(T_n(g_{r,s})\right)\) is performed element-wise and \(\mathbf{0}_{n_0}\) is a set of zeros of size \(n_0\) defined in \eqref{eq:n0}. An implementation is presented in Appendix~\ref{code:alleigs}.

\section{Examples and numerical experiments}
\label{sec:examples}

\subsection{Example 1: \(\{r,s\}=\{1,1\}\)}
As the first example we study the case \(\{r,s\}=\{1,1\}\), where we have a closed form formula for the eigenvalues for arbitrary \(n\); \(\lambda_j(T_n(g_{1,1}))=2\cos\left(\frac{j\pi}{n+1}\right)\), for \(j=1,\ldots,n\). 
In our algorithm, we have for \(\{r,s\}=\{1,1\}\) that \(\mathcal{M}=\left[\begin{array}{r}1\\2\end{array}\right]\) and \(\mathcal{P}=\left[\begin{array}{r}1\\0\end{array}\right]\).

For \(\beta_s=0\), that is \(n\) is even, we have
\begin{align*}
\quad B_{n/2}^{n,1,1}&=
\underbrace{\left[
    \begin{array}{rrrrrrr}
        1&\\
        1&1\\
        &\ddots&\ddots\\
        &&1&1&\\
        &&&1&1
    \end{array}
\right]}_{T_{n/2}^{\mathop\intercal}(e^{-\mathbf{i}\theta}c(\theta)^1)}
\underbrace{\left[
    \begin{array}{rrrrrrr}
        1&1\\
        &1&1\\
        &&\ddots&\ddots\\
        &&&1&1\\
        &&&&1
    \end{array}
\right]^1}_{T_{n/2}(e^{-\mathbf{i}\theta}c(\theta))^1}
=
\left[
    \begin{array}{rrrrrrr}
        \cca 1&1\\
        1&2&1\\
        &\ddots&\ddots&\ddots\\
        &&1&2&1\\
        &&&1&2
    \end{array}
\right].
\end{align*}
Indicated in blue is the top left corner, the element which makes this matrix not a pure Toeplitz matrix. For this matrix we have \(b_{1,1}(\theta)=2+2\cos(\theta)\) and \(\lambda_j(B_{n/2}^{n,1,1})=b_{1,1}(\theta_{j,n/2})\) where \(\theta_{j,n/2}=\frac{j2\pi}{n+1}, j=1,\ldots, n/2\); see for example~\cite{Bozzo_1995}. Thus, \(\lambda_+(T_n(g_{1,1}))=\sqrt{b_{1,1}(\theta_{j,n/2})}\), and the full spectrum of \(T_n(g_{1,1})\), as described in Section~\ref{sec:main:fullspectrum}, is given by \(\pm\lambda_+(T_n(g_{1,1}))\). Indeed, this is equivalent to the classical formula.

For \(\beta_s=1\), that is \(n\) is odd, we have one zero eigenvalue since \(n_0=1\) in \eqref{eq:n0} and
\begin{align*}
\quad B_{(n-1)/2}^{n,1,1}&=
\underbrace{\left[
    \begin{array}{rrrrrrr}
        2&1\\
        1&2&1\\
        &\ddots&\ddots&\ddots\\
        &&1&2&1\\
        &&&1&2
    \end{array}
\right]}_{T_{(n-1)/2}^{\mathop\intercal}(e^{-\mathbf{i}\theta}c(\theta)^2)}
\underbrace{\left[
    \begin{array}{rrrrrrr}
        1&1\\
        &1&1\\
        &&\ddots&\ddots\\
        &&&1&1\\
        &&&&1
    \end{array}
\right]^0}_{T_{(n-1)/2}(e^{-\mathbf{i}\theta}c(\theta))^0}=
\left[
    \begin{array}{rrrrrrr}
        2&1\\
        1&2&1\\
        &\ddots&\ddots&\ddots\\
        &&1&2&1\\
        &&&1&2
    \end{array}
\right].
\end{align*}
We notice that this matrix is pure Toeplitz, with the same symbol \(b_{1,1}(\theta)\) as above, but the eigenvalues \(\lambda_+(T_n(g_{1,1}))=\sqrt{b_{1,1}(\theta_{j,(n-1)/2})}\) are now given by the grid \(\theta_{j,(n-1)/2}=\frac{j2\pi}{n+1}, j=1,\ldots,(n-1)/2\). The eigenvalues of \(T_n(g_{1,1})\) are given by \(\pm\lambda_+(T_n(g_{1,1}))\cup 0\), and again this is equivalent to the classical formula.

\subsection{Example 2: \(\{r,s\}=\{1,2\}\)}
\label{examples:2}
\begin{figure}[H]
\centering
\includegraphics[width=0.425\textwidth]{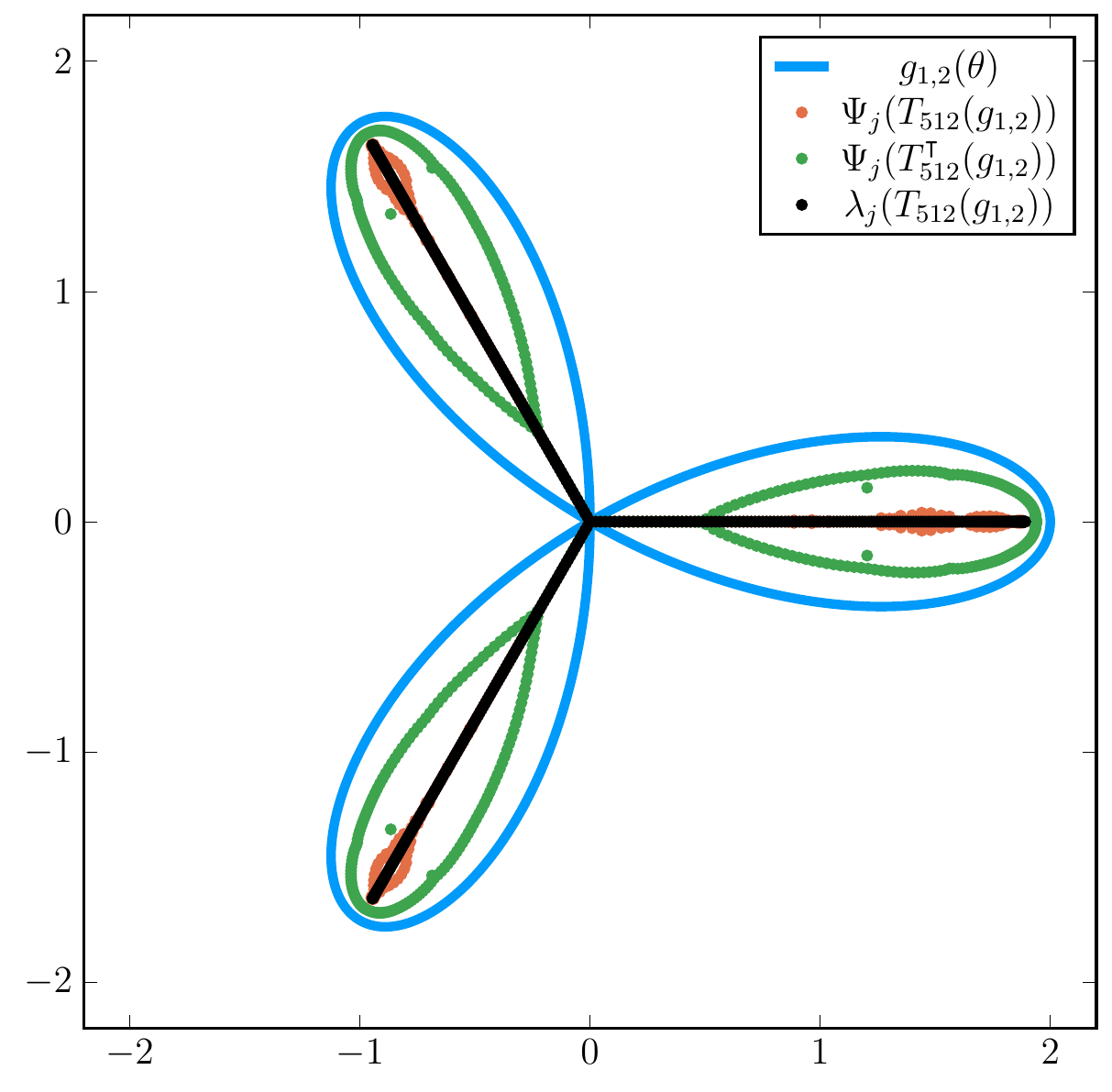}
\includegraphics[width=0.48\textwidth]{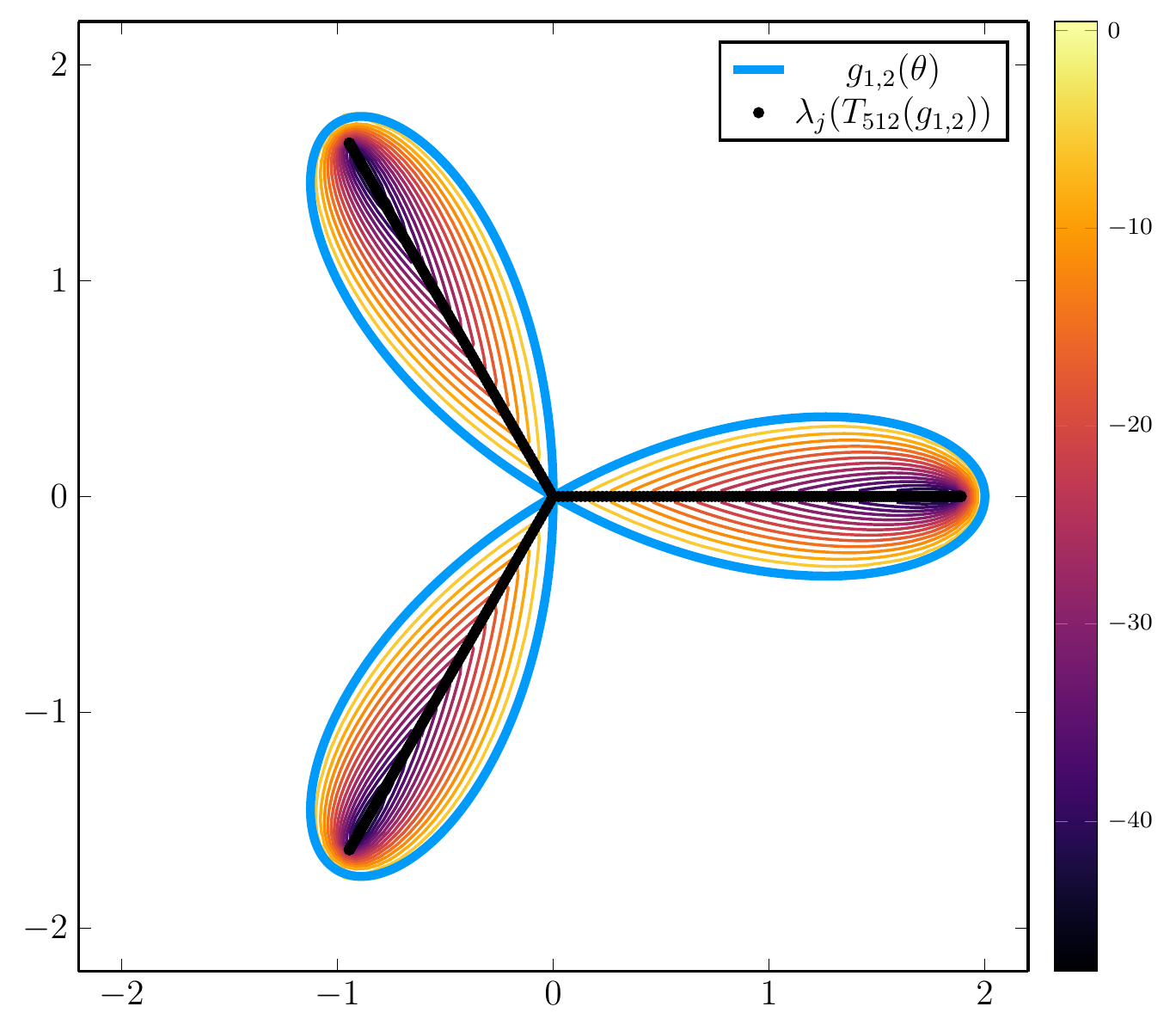}
\caption{Spectra and pseudospectra, \(\{n,r,s\}=\{512,1,2\}\): \textbf{Left:} Incorrect numerically computed eigenvalues \(\Psi_j(T_{512}(g_{1,2}))\) and \(\Psi_j(T_{512}^{\mathop\intercal}(g_{1,2}))\) using double precision and the  correctly computed eigenvalues \(\lambda_j(T_{512}(g_{1,2}))\) using 256-bit \textsc{BigFloat} in \textsc{Julia}. \textbf{Right:} The \(\epsilon\)-pseudospectrum of \(T_{512}(g_{1,2})\) computed using \textsc{Pseudospectra.jl}, where each contour is a different \(\log_{10}(\epsilon)\).}
    \label{fig:pseudospectra}
\end{figure}

As a second example we study the case \(\{r,s\}=\{1,2\}\). First we present the numerical difficulty of the nonsymmetric eigenvalue problem visually in Figure~\ref{fig:pseudospectra}. In the left panel we see the symbol \(g_{1,2}(\theta)\) (blue) and the true eigenvalues \(\lambda_j(T_{512}(g_{1,2}))\) (black) computed using high precision 256-bit \textsc{BigFloat} data type in \textsc{Julia}. Computing the spectrum using standard double precision, yields numerically unstable and incorrect results, presented as \(\Psi_j(T_{512}(g_{1,2}))\) (red) and \(\Psi_j(T_{512}^{\mathop\intercal}(g_{1,2}))\) (green). In the right panel we instead see the \(\epsilon\)-pseudospectrum~\cite{bottcher001,bpseudo,Bottcher2012-lf,Landau1975-fn,reichel921,trefethen051}, generated by the \textsc{Pseudospectra.jl}~\cite{pseudospectra.jl} package, a \textsc{Julia} implementation of the \textsc{MATLAB} toolbox \textsc{EigTool}~\cite{eigtool}.

In our algorithm, we have for \(\{r,s\}=\{1,2\}\) that \(\mathcal{M}=\left[\begin{array}{r}1\\2\\3\end{array}\right]\) and \(\mathcal{P}=\left[\begin{array}{r}2\\1\\0\end{array}\right]\).
Now, consider the three cases \(n=\{15, 16, 17\}\), that is,
\(n_\sigma=5\), and \(\beta_\sigma=\{0,1,2\}\).

\begin{figure}[H]
\centering
\includegraphics[width=0.5\textwidth]{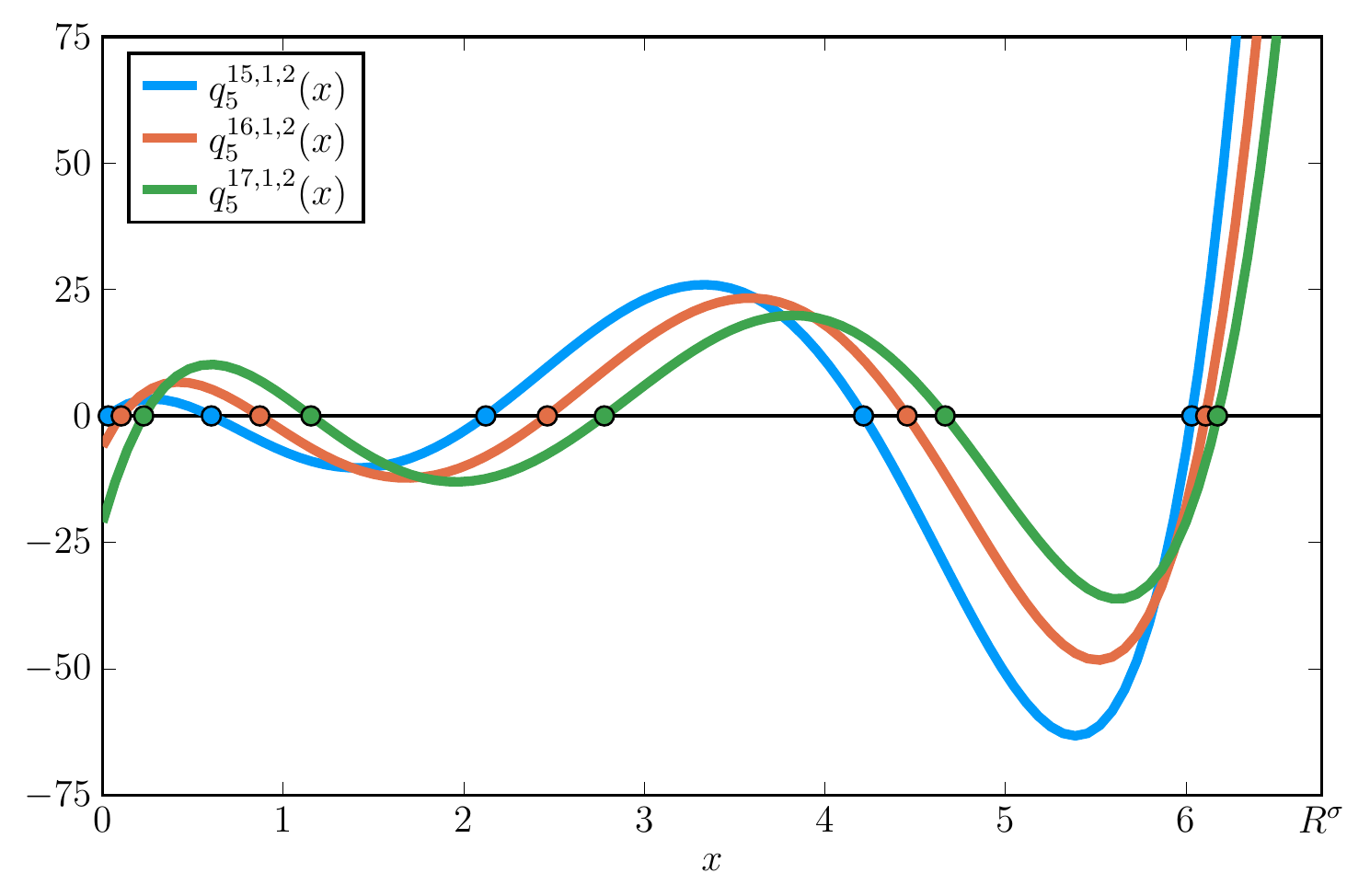}
\caption{The characteristic polynomials \(q_{n_\sigma}^{n,r,s}(x)\) for \(B_{n_\sigma}^{n,r,s}\) with \(n=\{15,16,17\}\), \(\{r,s\}=\{1,2\}\), and thus \(\beta_\sigma=0,1,2\). The roots are shown with circles. The positive real eigenvalues of \(T_n(g_{1,2})\) are in the interval \((0,R)\), where \(R^\sigma=\frac{\sigma^\sigma}{r^{r}s^{s}}=6.75\).}
\label{fig:charpoly12}
\end{figure}

Below we present the three generated \(B_{5}^{n,1,2}\), with the perturbations from the pure Toeplitz matrix \(T_{n_\sigma}(b_{1,2})\) in top left corners indicated in blue. 
\begin{align*}
   B_{5}^{15,1,2}=T_{5}^{\mathop\intercal}(e^{-\mathbf{i}\theta}c(\theta)^1)T_{5}(e^{-\mathbf{i}\theta}c(\theta))^2&=\left[
\begin{array}{ccccc}
    1 & 0 & 0 & 0 & 0 \\
    1 & 1 & 0 & 0 & 0 \\
    0 & 1 & 1 & 0 & 0 \\
    0 & 0 & 1 & 1 & 0 \\
    0 & 0 & 0 & 1 & 1 \\
\end{array}
\right]
\left[
\begin{array}{ccccc}
    1 & 1 & 0 & 0 & 0 \\
    0 & 1 & 1 & 0 & 0 \\
    0 & 0 & 1 & 1 & 0 \\
    0 & 0 & 0 & 1 & 1 \\
    0 & 0 & 0 & 0 & 1 \\
\end{array}
\right]^2=
    \left[
        \begin{array}{ccccc}
            \cca 1 & \cca 2 & 1 & 0 & 0 \\
            1 & 3 & 3 & 1 & 0 \\
            0 & 1 & 3 & 3 & 1 \\
            0 & 0 & 1 & 3 & 3 \\
            0 & 0 & 0 & 1 & 3 \\
        \end{array}
    \right],\\
B_{5}^{16,1,2}=T_{5}^{\mathop\intercal}(e^{-\mathbf{i}\theta}c(\theta)^2)T_{5}(e^{-\mathbf{i}\theta}c(\theta))^1&=\left[
        \begin{array}{ccccc}
        2 & 1 & 0 & 0 & 0 \\
        1 & 2 & 1 & 0 & 0 \\
        0 & 1 & 2 & 1 & 0 \\
        0 & 0 & 1 & 2 & 1 \\
        0 & 0 & 0 & 1 & 2 \\
        \end{array}
        \right]\left[
            \begin{array}{ccccc}
                1 & 1 & 0 & 0 & 0 \\
                0 & 1 & 1 & 0 & 0 \\
                0 & 0 & 1 & 1 & 0 \\
                0 & 0 & 0 & 1 & 1 \\
                0 & 0 & 0 & 0 & 1 \\
            \end{array}
            \right]^1
            =\left[
\begin{array}{ccccc}
    \cca 2 & 3 & 1 & 0 & 0 \\
    1 & 3 & 3 & 1 & 0 \\
    0 & 1 & 3 & 3 & 1 \\
    0 & 0 & 1 & 3 & 3 \\
    0 & 0 & 0 & 1 & 3 \\
\end{array}
\right],\\
B_{5}^{17,1,2}=T_{5}^{\mathop\intercal}(e^{-\mathbf{i}\theta}c(\theta)^3)T_{5}(e^{-\mathbf{i}\theta}c(\theta))^0&=\left[
\begin{array}{ccccc}
    3 & 3 & 1 & 0 & 0 \\
    1 & 3 & 3 & 1 & 0 \\
    0 & 1 & 3 & 3 & 1 \\
    0 & 0 & 1 & 3 & 3 \\
    0 & 0 & 0 & 1 & 3 \\
\end{array}
\right]\left[
    \begin{array}{ccccc}
        1 & 1 & 0 & 0 & 0 \\
        0 & 1 & 1 & 0 & 0 \\
        0 & 0 & 1 & 1 & 0 \\
        0 & 0 & 0 & 1 & 1 \\
        0 & 0 & 0 & 0 & 1 \\
    \end{array}
    \right]^0=\left[
\begin{array}{ccccc}
    3 & 3 & 1 & 0 & 0 \\
    1 & 3 & 3 & 1 & 0 \\
    0 & 1 & 3 & 3 & 1 \\
    0 & 0 & 1 & 3 & 3 \\
    0 & 0 & 0 & 1 & 3 \\
\end{array}
\right].
\end{align*}
The number of zero eigenvalues of \(T_n(g_{1,2})\) are \(n_0=\{0,1,2\}\) for \(n=\{15,16,17\}\) respectively. 
Furthermore, we can visualise how the eigenvalues change for different \(\beta_\sigma\) by looking at the characteristic polynomials \(q_{n_\sigma}^{n,r,s}(x)\) for \(B_{n_\sigma}^{n,r,s}\), that is \(n=\{15,16,17\}\) for \(\beta_\sigma=0,1,2\), are
\begin{align*}
    q_5^{15,1,2}(x)&=x^5 - 13x^4 + 55x^3 - 84x^2 + 35x - 1\\
    q_5^{16,1,2}(x)&=x^5 - 14x^4 + 66x^3 - 120x^2 + 70x - 6\\
    q_5^{17,1,2}(x)&=x^5 - 15x^4 + 78x^3 - 165x^2 + 126x - 21
\end{align*}
and the roots (which are the eigenvalues \(\lambda_+(T_n(g_{1,2}))^\sigma\)) are shown in Figure~\ref{fig:charpoly12}. The smallest eigenvalue is given by \(q_5^{15,1,2}(x)\), that is, \(\beta_\sigma=0\), resulting in a larger condition number \(\kappa(\cdot)=\frac{\sigma{\mathrm{max}}}{\sigma{\mathrm{min}}}\) than for \(\beta_\sigma>0\); \(\kappa=\{242.8,87.4,38.3\}\). Also, we notice this pattern in Section~\ref{sec:experiments:numerical}, that the error and \(\kappa(B_{n_\sigma}^{n,r,s})\) is declining as \(\beta_{\sigma}\) increases; See left panel in Figure~\ref{fig:Error}.
The characteristic polynomials for \(T_n(g_{r,s})\) are \(p_n^{r,s}(x)\). As expected, we have
\begin{align*}
    p_{15}^{1,2}(x)&=q_5^{15,1,2}(x^3),\\
    p_{16}^{1,2}(x)&=xq_5^{16,1,2}(x^3),\\
    p_{17}^{1,2}(x)&=x^2q_5^{17,1,2}(x^3),
\end{align*}
in accordance with Conjecture~\ref{conj:1}.

\begin{rmrk}
We note that \(R\) defined in  \cite[Section 7]{schmidt601} to the power of \(\sigma\) is in this case defined by
\begin{align*}
R^\sigma=\frac{\sigma^\sigma}{r^{r}s^{s}}=\frac{3^3}{1^1 2^2}=\frac{27}{4}=6.75,
\end{align*} 
so  the upper bound of the real eigenvalues of \(T_n(g_{1,2})\) is \(R=6.75^{\frac{1}{3}}\).
For the general case of any \{r,s\}, and \(\gamma=1\), we have that \(\mathfrak{b}_{r,s}\) in \eqref{eq:shapirob},
\begin{align*}
    \lim_{\theta\to 0}\mathfrak{b}_{r,s}(\theta)=\lim_{\theta\to 0}\frac{\sin^\sigma(\theta)}{\sin^r(\frac{r}{\sigma}\theta)\sin^s(\frac{s}{\sigma}\theta)}
    =\lim_{\theta\to 0}\frac{\left(\frac{\sigma}{\sigma}\theta\right)^\sigma}{\left(\frac{r}{\sigma}\theta\right)^r\left(\frac{s}{\sigma}\theta\right)^s}
    =\frac{\sigma^\sigma}{r^rs^s}
    =R^\sigma,
\end{align*}
since \(\mathfrak{b}_{r,s}\) is even, and monotone in \([0,\pi]\), and has its maximum at \(\theta=0\).
\end{rmrk}

\subsection{Example 3: \(\{r,s\}=\{2,4\}\)}
 For \(\{r,s\}=\{2,4\}\) we have that \(\sigma=6\), \(\gamma=2\), \(\{r_\gamma,s_\gamma\}=\{1,2\}\), and \(\sigma_\gamma=3\). Then for a given \(n\), the eigenvalues of \(T_n(g_{2,4})\) will be computed by using the union of the eigenvalues of \(B_{(n_\gamma)_\sigma}^{n_\gamma,1,2}\) repeating them \(\gamma - \beta_\gamma\) times and \(B_{(n_\gamma+1)_\sigma}^{n_\gamma+1,1,2}\) repeated \(\beta_\gamma\) times. Take the cases \(n=\{12, 13, 14, 15,16,17\}\), where the relevant variables have been computed in Table~\ref{table:example3}.
\begin{table}[H]
    \centering
    \caption{Variables for Example 3.}\label{table:example3}
    \begin{tabular}{c|cc|cc|ccc|c}
    \toprule
    \(n\) & \(\beta_\sigma\) & \(n_\sigma\) & \(\beta_\gamma\) & \(n_\gamma\) & \((\beta_\gamma)_\sigma\) & \((n_\gamma)_\sigma\) & \((n_\gamma + 1)_\sigma\) & \(n_0\)\\ \midrule
    12 & 0 & 2 & 0 & 6 & 0 & 2 & 2 & 0 \\
    13 & 1 & 2 & 1 & 6 & 0 & 2 & 2 & 1 \\
    14 & 2 & 2 & 0 & 7 & 1 & 2 & 2 & 2 \\
    15 & 3 & 2 & 1 & 7 & 1 & 2 & 2 & 3 \\
    16 & 4 & 2 & 0 & 8 & 2 & 2 & 3 & 4 \\
    17 & 5 & 2 & 1 & 8 & 2 & 2 & 3 & 2 \\ \bottomrule
    \end{tabular}
\end{table}

\noindent Thus, the necessary matrices to compute the eigenvalues are,
\begin{align}
B_2^{6, 1, 2}&=\left[
    \begin{array}{cc}
    1 & 2 \\
    1 & 3 \\
    \end{array}
    \right],\qquad
B_2^{7, 1, 2}=\left[
    \begin{array}{cc}
    2 & 3 \\
    1 & 3 \\
    \end{array}
    \right],\qquad
B_2^{8, 1, 2}=\left[
    \begin{array}{cc}
    3 & 3 \\
    1 & 3 \\
    \end{array}
    \right],\qquad
B_3^{9, 1, 2}=\left[
    \begin{array}{ccc}
    1 & 2 & 1 \\
    1 & 3 & 3 \\
    0 & 1 & 3 \\
    \end{array}
    \right],\label{eq:example3:Bnmats}
\end{align}
so, we compute the positive real eigenvalues using
\begin{align*}
    \lambda_+(T_{12}(g_{2,4}))^3 &= \bigcup_{k=1}^{2}\left( \bigcup_{j=1}^{2}\lambda_j\left(B_2^{6, 1, 2}\right)\right),\\
    \lambda_+(T_{13}(g_{2,4}))^3 &= \left(\bigcup_{j=1}^{2}\lambda_j\left(B_2^{6, 1, 2}\right)\right) \bigcup \left(\bigcup_{j=1}^{2}\lambda_j\left(B_2^{7, 1, 2}\right)\right),\\
    \lambda_+(T_{14}(g_{2,4}))^3 &= \bigcup_{k=1}^{2}\left(\bigcup_{j=1}^{2} \lambda_j\left(B_2^{7, 1, 2}\right)\right),\\
    \lambda_+(T_{15}(g_{2,4}))^3 &= \left(\bigcup_{j=1}^{2}\lambda_j\left(B_2^{7, 1, 2}\right)\right) \bigcup \left(\bigcup_{j=1}^{2}\lambda_j\left(B_2^{8, 1, 2}\right)\right),\\
    \lambda_+(T_{16}(g_{2,4}))^3 &= \bigcup_{k=1}^{2}\left(\bigcup_{j=1}^{2} \lambda_j\left(B_2^{8, 1, 2}\right)\right),\\
    \lambda_+(T_{17}(g_{2,4}))^3 &= \left(\bigcup_{j=1}^{2}\lambda_j\left(B_2^{8, 1, 2}\right)\right) \bigcup \left(\bigcup_{j=1}^{3}\lambda_j\left(B_3^{9, 1, 2}\right)\right),
\end{align*}
where the eigenvalues of the matrices \eqref{eq:example3:Bnmats} are either computed analytically or numerically.
We see that in general we need to use a union of the eigenvalues of a total of \(\gamma\) matrices (where it is either one or two different ones) to compute the positive real eigenvalues of \(T_{n}(g_{r,s})\).

In Table~\ref{tbl:exampl3:charpoly} is shown the characteristic polynomial of \(T_{n}(g_{2,4})\) for \(n=\{12, 13, 14, 15,16,17\}\) as well as the matrices constructed in \eqref{eq:example3:Bnmats}.

\begin{table}[H]
    \centering
    \caption{Characteristic polynomials of matrices in Example 3.}
    \label{tbl:exampl3:charpoly}
\begin{tabular}{rl}
\toprule 
Matrix&Characteristic polynomial\\
\midrule
\(T_{12}(g_{2,4})\)&\(p_{12}^{2,4}(x)=x^{12} - 8x^9 + 18x^6 - 8x^3 + 1\)\\
\(T_{13}(g_{2,4})\)&\(p_{13}^{2,4}(x)=x(x^{12} - 9x^{9} + 24x^6 - 17x^3 + 3)\)\\
\(T_{14}(g_{2,4})\)&\(p_{14}^{2,4}(x)=x^2(x^{12} - 10x^{9} + 31x^6 - 30x^3 + 9)\)\\
\(T_{15}(g_{2,4})\)&\(p_{15}^{2,4}(x)=x^3(x^{12} - 11x^{9} + 39x^6 - 48x^3 + 18)\)\\
\(T_{16}(g_{2,4})\)&\(p_{16}^{2,4}(x)=x^4(x^{12} - 12x^{9} + 48x^{6} - 72x^3 + 36)\)\\
\(T_{17}(g_{2,4})\)&\(p_{17}^{2,4}(x)=x^2(x^{15} - 13x^{12} + 58x^{9} - 103x^6 + 66x^3 - 6)\)\\
\(B_2^{6, 1, 2}\)&\(q_2^{6, 1, 2}(x)=x^2 - 4x + 1\)\\
\(B_2^{7, 1, 2}\)&\(q_2^{7, 1, 2}(x)=x^2 - 5x + 3\)\\
\(B_2^{8, 1, 2}\)&\(q_2^{8, 1, 2}(x)=x^2 - 6x + 6\)\\
\(B_3^{9, 1, 2}\)&\(q_3^{9, 1, 2}(x)=x^3 - 7x^2 + 10x - 1\)\\
\bottomrule
\end{tabular}
\end{table}
We have from Table~\ref{tbl:exampl3:charpoly} and the number of zeros \(n_0\) in Table~\ref{table:example3},
\begin{align*}
    p_{12}^{2,4}(x)&=q_2^{6,1,2}(x^3)q_2^{6,1,2}(x^3),\\
    p_{13}^{2,4}(x)&=xq_2^{6,1,2}(x^3)q_2^{7,1,2}(x^3),\\
    p_{14}^{2,4}(x)&=x^2q_2^{7,1,2}(x^3)q_2^{7,1,2}(x^3),\\
    p_{15}^{2,4}(x)&=x^3q_2^{7,1,2}(x^3)q_2^{8,1,2}(x^3),\\
    p_{16}^{2,4}(x)&=x^4q_2^{8,1,2}(x^3)q_2^{8,1,2}(x^3),\\
    p_{17}^{2,4}(x)&=x^2q_2^{8,1,2}(x^3)q_3^{9,1,2}(x^3),\\
\end{align*}
so indeed, Conjecture~\ref{conj:1} is verified also in this case.
Furthermore, we note that, for example,
\begin{align*}
T_{12}(g_{2,4})^3=
\left[
\begin{array}{cccccccccccc}
1 & 0 & 0 & 0 & 0 & 0 & 2 & 0 & 0 & 0 & 0 & 0 \\
0 & 1 & 0 & 0 & 0 & 0 & 0 & 2 & 0 & 0 & 0 & 0 \\
0 & 0 & 2 & 0 & 0 & 0 & 0 & 0 & 3 & 0 & 0 & 0 \\
0 & 0 & 0 & 2 & 0 & 0 & 0 & 0 & 0 & 3 & 0 & 0 \\
0 & 0 & 0 & 0 & 3 & 0 & 0 & 0 & 0 & 0 & 2 & 0 \\
0 & 0 & 0 & 0 & 0 & 3 & 0 & 0 & 0 & 0 & 0 & 2 \\
1 & 0 & 0 & 0 & 0 & 0 & 3 & 0 & 0 & 0 & 0 & 0 \\
0 & 1 & 0 & 0 & 0 & 0 & 0 & 3 & 0 & 0 & 0 & 0 \\
0 & 0 & 1 & 0 & 0 & 0 & 0 & 0 & 2 & 0 & 0 & 0 \\
0 & 0 & 0 & 1 & 0 & 0 & 0 & 0 & 0 & 2 & 0 & 0 \\
0 & 0 & 0 & 0 & 1 & 0 & 0 & 0 & 0 & 0 & 1 & 0 \\
0 & 0 & 0 & 0 & 0 & 1 & 0 & 0 & 0 & 0 & 0 & 1 \\
\end{array}
\right]=
P^{-1}\left[
\begin{array}{cccccccccccc}
    1 & 2 &  &  &  &  &  &  &  &  &  &  \\
    1 & 3 &  &  &  &  &  &  &  &  &  &  \\
     &  & 1 & 2 &  &  &  &  &  &  &  &  \\
     &  & 1 & 3 &  &  &  &  &  &  &  &  \\
     &  &  &  & 1 & 1 &  &  &  &  &  &  \\
     &  &  &  & 2 & 3 &  &  &  &  &  &  \\
     &  &  &  &  &  & 1 & 1 &  &  &  &  \\
     &  &  &  &  &  & 2 & 3 &  &  &  &  \\
     &  &  &  &  &  &  &  & 2 & 1 &  &  \\
     &  &  &  &  &  &  &  & 3 & 2 &  &  \\
     &  &  &  &  &  &  &  &  &  & 2 & 1 \\
     &  &  &  &  &  &  &  &  &  & 3 & 2 \\
\end{array}
\right]P,
\end{align*}
where \(P\) is a similarity transformation, and two blocks are \(B_2^{6, 1, 2}\), two blocks are \((B_2^{6, 1, 2})^{\mathop\intercal}\), and two blocks are similar to \(B_2^{6, 1, 2}\). Also, to show an example with a zero eigenvalue we have 
\begin{align*}
    T_{13}(g_{2,4})^3=
    \left[
\begin{array}{ccccccccccccc}
1 & 0 & 0 & 0 & 0 & 0 & 2 & 0 & 0 & 0 & 0 & 0 & 1 \\
0 & 1 & 0 & 0 & 0 & 0 & 0 & 2 & 0 & 0 & 0 & 0 & 0 \\
0 & 0 & 2 & 0 & 0 & 0 & 0 & 0 & 3 & 0 & 0 & 0 & 0 \\
0 & 0 & 0 & 2 & 0 & 0 & 0 & 0 & 0 & 3 & 0 & 0 & 0 \\
0 & 0 & 0 & 0 & 3 & 0 & 0 & 0 & 0 & 0 & 3 & 0 & 0 \\
0 & 0 & 0 & 0 & 0 & 3 & 0 & 0 & 0 & 0 & 0 & 2 & 0 \\
1 & 0 & 0 & 0 & 0 & 0 & 3 & 0 & 0 & 0 & 0 & 0 & 2 \\
0 & 1 & 0 & 0 & 0 & 0 & 0 & 3 & 0 & 0 & 0 & 0 & 0 \\
0 & 0 & 1 & 0 & 0 & 0 & 0 & 0 & 3 & 0 & 0 & 0 & 0 \\
0 & 0 & 0 & 1 & 0 & 0 & 0 & 0 & 0 & 2 & 0 & 0 & 0 \\
0 & 0 & 0 & 0 & 1 & 0 & 0 & 0 & 0 & 0 & 2 & 0 & 0 \\
0 & 0 & 0 & 0 & 0 & 1 & 0 & 0 & 0 & 0 & 0 & 1 & 0 \\
0 & 0 & 0 & 0 & 0 & 0 & 1 & 0 & 0 & 0 & 0 & 0 & 1 \\
\end{array}
\right]=
    P^{-1}\left[
    \begin{array}{ccccccccccccc}
        1 & 2 &  &  &  &  &  &  &  &  &  &  &  \\
        1 & 3 &  &  &  &  &  &  &  &  &  &  &  \\
         &  & 1 & 1 &  &  &  &  &  &  &  &  &  \\
         &  & 2 & 3 &  &  &  &  &  &  &  &  &  \\
         &  &  &  & 2 & 1 &  &  &  &  &  &  &  \\
         &  &  &  & 3 & 2 &  &  &  &  &  &  &  \\
         &  &  &  &  &  & 2 & 3 &  &  &  &  &  \\
         &  &  &  &  &  & 1 & 3 &  &  &  &  &  \\
         &  &  &  &  &  &  &  & 2 & 1 &  &  &  \\
         &  &  &  &  &  &  &  & 3 & 3 &  &  &  \\
         &  &  &  &  &  &  &  &  &  & 1 & 2 & 1 \\
         &  &  &  &  &  &  &  &  &  & 1 & 3 & 2 \\
         &  &  &  &  &  &  &  &  &  & 0 & 1 & 1 \\
    \end{array}
    \right]P,
    \end{align*}
where one blocks is \(B_2^{6, 1, 2}\), one block is \((B_2^{6, 1, 2})^{\mathop\intercal}\), and one block is similar to \(B_2^{6, 1, 2}\). Then, one blocks is \(B_2^{7, 1, 2}\), one block is \((B_2^{7, 1, 2})^{\mathop\intercal}\), and lastly a \(3\times 3\) block with same eigenvalues as \(B_2^{7, 1, 2}\) plus a zero eigenvalue.

\subsection{Example 4: \(\{r,s\} = \{3,5\}\)}
\label{sec:example:CornerAndPerms}
In this fourth example, we wish to show how the matrices \(\mathcal{M}\) and \(\mathcal{P}\) are permuted as well as how the erroneous corner elements are handled when \(\beta_\sigma > s\). Consider the cases of matrices generated for \(\{n,r,s\}\) being \(\{86,3,5\}\) and \(\{87,3,5\}\); that is \(n_\sigma=10\) and \(\beta_\sigma=6,7\). 
As we have \(r = 3 > 1\), we will need to permute the intiial \(\mathcal{M}\) and \(\mathcal{P}\) matrices to get the correct \(\mathcal{M}\) and \(\mathcal{P}\). The initial constructed Toeplitz matrices are
\begin{align*}
    \mathcal{M}=
    \left[
    \begin{array}{rrr}
    1 & 1 & 1 \\
    2 & 1 & 1 \\
    2 & 2 & 1 \\
    2 & 2 & 2 \\
    3 & 2 & 2 \\
    3 & 3 & 2 \\
    3 & 3 & 3 \\
    4 & 3 & 3 \\
    \end{array}
    \right],
    \qquad
    \mathcal{P}=
    \left[
    \begin{array}{rrr}
    1 & 2 & 2 \\
    1 & 1 & 2 \\
    1 & 1 & 1 \\
    0 & 1 & 1 \\
    0 & 0 & 1 \\
    0 & 0 & 0 \\
    -1 & 0 & 0 \\
    -1 & -1 & 0 \\
    \end{array}
    \right],
    \end{align*}
which we will then permute using \(p_\mathrm{perm}\). First we compute \(\tau = \bmod(s,r) = 2\), then we construct the ordered set \(\mathcal{S}_\tau = \left(\tau, 2\tau, \dots, r\tau\right) = \left(2, 4, 6\right)\) which we then use to compute \(p_\mathrm{perm} = \bmod(\mathcal{S}_\tau, r) = (2, 1, 0)\); however the last element should be equal to \(r=3\), and hence the correct permutation is \(p_\mathrm{perm} = (2, 1, 3)\).   The permuted, and correct, matrices \(\mathcal{M}\) and \(\mathcal{P}\) for \(\{r,s\}=\{3,5\}\) are then
\begin{align*}
\mathcal{M}=
\left[
\begin{array}{rrr}
1 & 1 & 1 \\
1 & 2 & 1 \\
2 & 2 & 1 \\
2 & 2 & 2 \\
2 & 3 & 2 \\
3 & 3 & 2 \\
3 & 3 & 3 \\
3 & 4 & 3 \\
\end{array}
\right],
\qquad
\mathcal{P}=
\left[
\begin{array}{rrr}
2 & 1 & 2 \\
1 & 1 & 2 \\
1 & 1 & 1 \\
1 & 0 & 1 \\
0 & 0 & 1 \\
0 & 0 & 0 \\
0 & -1 & 0 \\
-1 & -1 & 0 \\
\end{array}
\right].
\end{align*}
As we are in the case \(\beta_\sigma=6,7\), we have that \(\beta_\sigma > s\), and thus we will have erroneous corner elements. In this setting, it is the last two rows of \(\mathcal{M}\) and \(\mathcal{P}\) that define the construction of \(B_{n_\sigma}^{n,r,s}\), and includes negative values in \(\mathcal{P}\).
To get the correct \( B_{n_\sigma}^{n,r,s}\) we replace the red elements on the right top corner with zeros.
\begin{align*}
    B_{10}^{86,3,5}-R_{10}^{86,3,5}&=
    T_{n_\sigma}^{\mathop\intercal}(e^{-\mathbf{i}\theta}c(\theta)^3)
    T_{n_\sigma}^{\mathop\intercal}(e^{-\mathbf{i}\theta}c(\theta)^3)
    (T_{n_\sigma}(e^{-\mathbf{i}\theta}c(\theta)))^{-1}
    T_{n_\sigma}^{\mathop\intercal}(e^{-\mathbf{i}\theta}c(\theta)^3)\\
    &=
    \left[
        \begin{array}{cccccccccc}
            43 & 65 & 55 & 28 & 8 & 1 & 0 & 0 & 0 & \ccb -3 \\
            27 & 56 & 70 & 56 & 28 & 8 & 1 & 0 & 0 & \ccb -1 \\
            8 & 28 & 56 & 70 & 56 & 28 & 8 & 1 & 0 & 0 \\
            1 & 8 & 28 & 56 & 70 & 56 & 28 & 8 & 1 & 0 \\
            0 & 1 & 8 & 28 & 56 & 70 & 56 & 28 & 8 & 1 \\
            0 & 0 & 1 & 8 & 28 & 56 & 70 & 56 & 28 & 8 \\
            0 & 0 & 0 & 1 & 8 & 28 & 56 & 70 & 56 & 28 \\
            0 & 0 & 0 & 0 & 1 & 8 & 28 & 56 & 70 & 55 \\
            0 & 0 & 0 & 0 & 0 & 1 & 8 & 28 & 55 & 62 \\
            0 & 0 & 0 & 0 & 0 & 0 & 1 & 8 & 25 & 37 \\
        \end{array}
        \right],\\
        B_{10}^{87,3,5}-R_{10}^{87,3,5}&=  
        T_{n_\sigma}^{\mathop\intercal}(e^{-\mathbf{i}\theta}c(\theta)^3)
        (T_{n_\sigma}^{\mathop\intercal}(e^{-\mathbf{i}\theta}c(\theta)))^{-1}
        T_{n_\sigma}^{\mathop\intercal}(e^{-\mathbf{i}\theta}c(\theta)^4)
        (T_{n_\sigma}^{\mathop\intercal}(e^{-\mathbf{i}\theta}c(\theta)))^{-1}
        T_{n_\sigma}^{\mathop\intercal}(e^{-\mathbf{i}\theta}c(\theta)^3)\\
        &=
        \left[
\begin{array}{cccccccccc}
    43 & 65 & 55 & 28 & 8 & 1 & 0 & 0 & \ccb -1 & \ccb -6 \\
    27 & 56 & 70 & 56 & 28 & 8 & 1 & 0 & 0 & \ccb -1 \\
    8 & 28 & 56 & 70 & 56 & 28 & 8 & 1 & 0 & 0 \\
    1 & 8 & 28 & 56 & 70 & 56 & 28 & 8 & 1 & 0 \\
    0 & 1 & 8 & 28 & 56 & 70 & 56 & 28 & 8 & 1 \\
    0 & 0 & 1 & 8 & 28 & 56 & 70 & 56 & 28 & 8 \\
    0 & 0 & 0 & 1 & 8 & 28 & 56 & 70 & 56 & 28 \\
    0 & 0 & 0 & 0 & 1 & 8 & 28 & 56 & 70 & 55 \\
    0 & 0 & 0 & 0 & 0 & 1 & 8 & 28 & 56 & 65 \\
    0 & 0 & 0 & 0 & 0 & 0 & 1 & 8 & 27 & 43 \\
\end{array}
\right].
\end{align*}
We have thus identified the low rank matrix \(R_{n_\sigma}^{n,r,s}\) (for all matrices \(B_{n_\sigma}^{n,r,s}\) large enough with the same \(\beta_\sigma\)) and can use this  to construct the correct \(B_{n_\sigma}^{n,r,s}\) for smaller matrices; see the four examples below.
\begin{align*}
    B_{4}^{38,3,5}-R_{4}^{38,3,5}&=\left[
        \begin{array}{cccc}
            43 & 65 & 55 & \ccb 25 \\
            27 & 56 & 70 & \ccb 54 \\
            8 & 28 & 55 & 62 \\
            1 & 8 & 25 & 37 \\
        \end{array}
        \right]=
    \underbrace{\left[
\begin{array}{cccc}
    43 & 65 & 55 & 28 \\
    27 & 56 & 70 & 55 \\
    8 & 28 & 55 & 62 \\
    1 & 8 & 25 & 37 \\
\end{array}
\right]}_{B_{4}^{38,3,5}}-\underbrace{\left[
    \begin{array}{cccc}
    0 & 0 & 0 & \ccb 3 \\
    0 & 0 & 0 & \ccb 1 \\
    0 & 0 & 0 & 0 \\
    0 & 0 & 0 & 0 \\
    \end{array}
    \right]}_{R_{4}^{38,3,5}}\\
    B_{4}^{39,3,5}-R_{4}^{39,3,5}&=\left[
        \begin{array}{cccc}
        43 & 65 & \ccb 54 & \ccb 22 \\
        27 & 56 & 70 & \ccb 54 \\
        8 & 28 & 56 & 65 \\
        1 & 8 & 27 & 43 \\
        \end{array}
        \right]=
    \underbrace{\left[
\begin{array}{cccc}
43 & 65 & 55 & 28 \\
27 & 56 & 70 & 55 \\
8 & 28 & 56 & 65 \\
1 & 8 & 27 & 43 \\
\end{array}
\right]}_{B_{4}^{39,3,5}}-\underbrace{\left[
    \begin{array}{cccc}
    0 & 0 & \ccb 1 & \ccb 6 \\
    0 & 0 & 0 & \ccb 1 \\
    0 & 0 & 0 & 0 \\
    0 & 0 & 0 & 0 \\
    \end{array}
    \right]}_{R_{4}^{39,3,5}}
\end{align*}
Note that we add \(R_{n_\sigma}^{n,r,s}\) from the originally constructed \(B_{n_\sigma}^{n,r,s}\) to get the true \(B_{n_\sigma}^{n,r,s}\), since  \(n_\sigma=4\) is even. 
\begin{align*}
B_{5}^{46,3,5}+R_{5}^{46,3,5}&=
\left[
\begin{array}{ccccc}
43 & 65 & 55 & 28 & \ccb 11 \\
27 & 56 & 70 & 56 & \ccb 29 \\
8 & 28 & 56 & 70 & 55 \\
1 & 8 & 28 & 55 & 62 \\
0 & 1 & 8 & 25 & 37 \\
\end{array}
\right]=
\underbrace{\left[
\begin{array}{ccccc}
43 & 65 & 55 & 28 & 8 \\
27 & 56 & 70 & 56 & 28 \\
8 & 28 & 56 & 70 & 55 \\
1 & 8 & 28 & 55 & 62 \\
0 & 1 & 8 & 25 & 37 \\
\end{array}
\right]}_{B_{5}^{46,3,5}}+
\underbrace{\left[
    \begin{array}{ccccc}
    0 & 0 & 0 & 0 & \ccb 3 \\
    0 & 0 & 0 & 0 & \ccb 1 \\
    0 & 0 & 0 & 0 & 0 \\
    0 & 0 & 0 & 0 & 0 \\
    0 & 0 & 0 & 0 & 0 \\
    \end{array}
    \right]}_{R_{5}^{46,3,5}}\\
    B_{5}^{47,3,5}+R_{5}^{47,3,5}&=
    \left[
    \begin{array}{ccccc}
        43 & 65 & 55 & \ccb 29 & \ccb 14 \\
        27 & 56 & 70 & 56 & \ccb 29 \\
        8 & 28 & 56 & 70 & 55 \\
        1 & 8 & 28 & 56 & 65 \\
        0 & 1 & 8 & 27 & 43 \\
    \end{array}
    \right]=
    \underbrace{\left[
    \begin{array}{ccccc}
        43 & 65 & 55 & 28 & 8 \\
        27 & 56 & 70 & 56 & 28 \\
        8 & 28 & 56 & 70 & 55 \\
        1 & 8 & 28 & 56 & 65 \\
        0 & 1 & 8 & 27 & 43 \\
    \end{array}
    \right]}_{B_{5}^{47,3,5}}+
    \underbrace{\left[
        \begin{array}{ccccc}
        0 & 0 & 0 & \ccb 1 & \ccb 6 \\
        0 & 0 & 0 & 0 & \ccb 1 \\
        0 & 0 & 0 & 0 & 0 \\
        0 & 0 & 0 & 0 & 0 \\
        0 & 0 & 0 & 0 & 0 \\
        \end{array}
        \right]}_{R_{5}^{47,3,5}}
\end{align*}
Note that we subtract \(R_{n_\sigma}^{n,r,s}\) from the originally constructed \(B_{n_\sigma}^{n,r,s}\) to get the true \(B_{n_\sigma}^{n,r,s}\), since  \(n_\sigma=5\) is odd. 

\subsection{Numerical experiments}
\label{sec:experiments:numerical}
\begin{figure}[H]
    \centering
    \includegraphics[width=.48\textwidth]{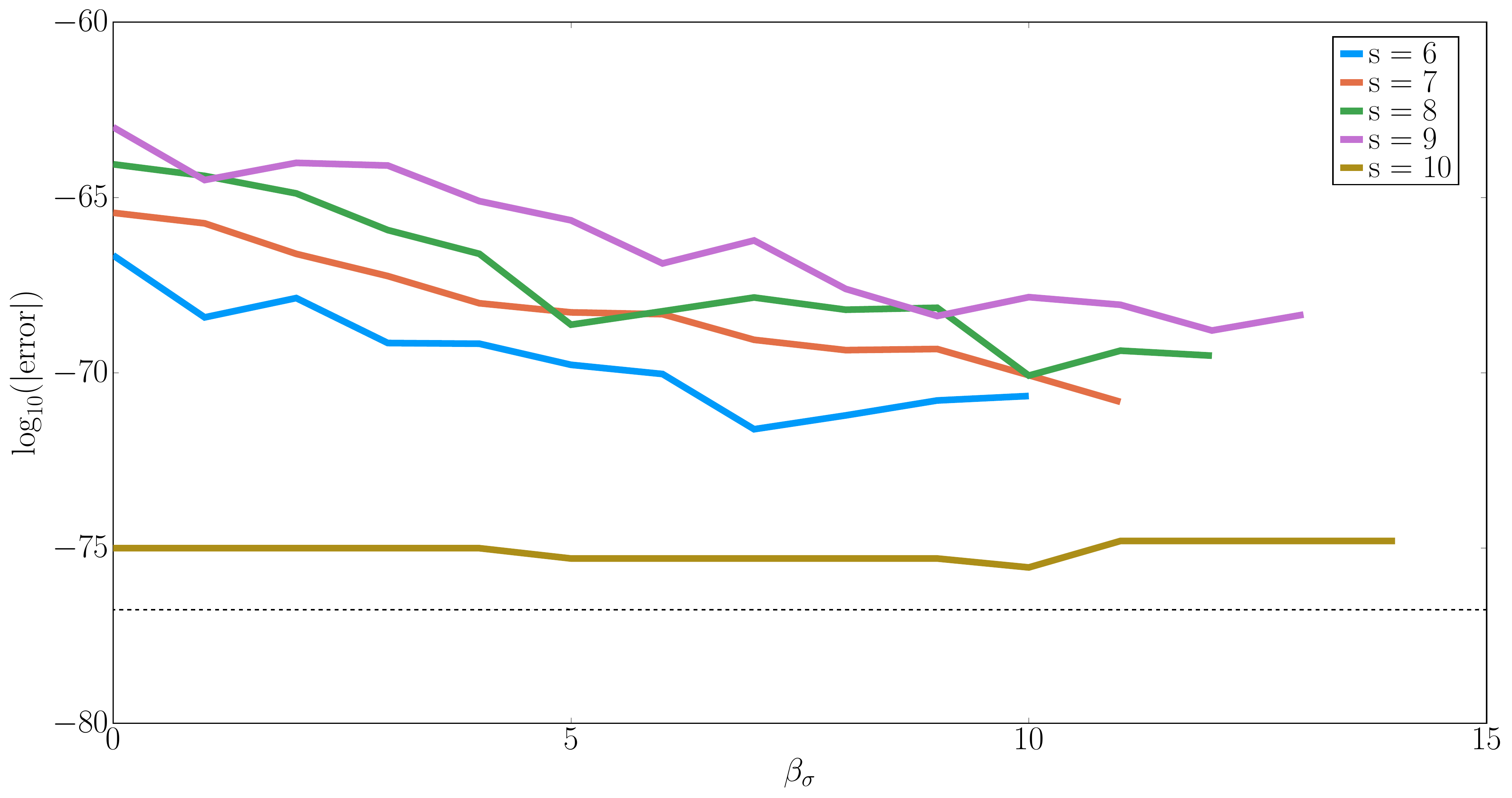}
    \includegraphics[width=.49\textwidth]{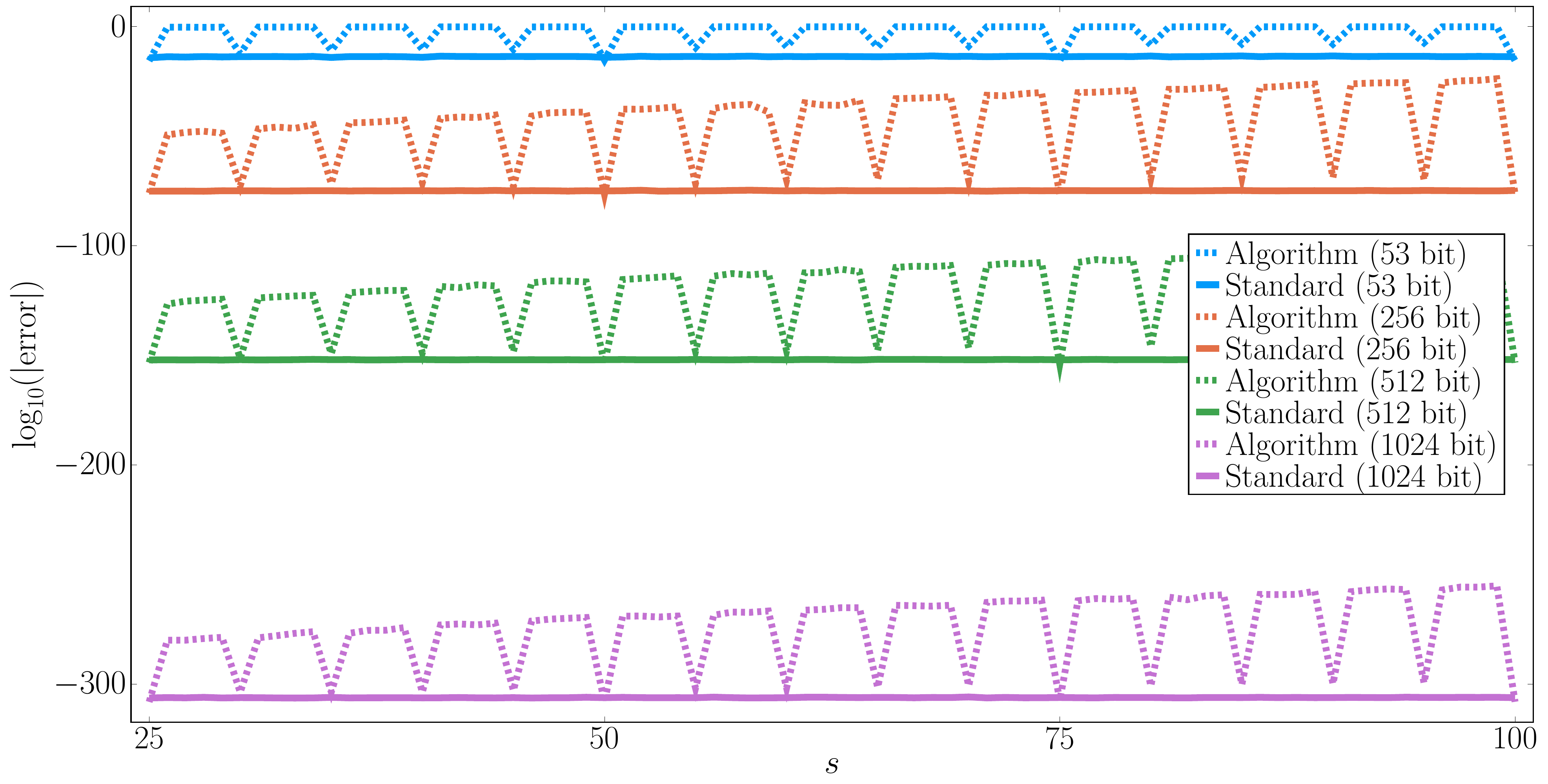}
    \caption{\textbf{Left:} Errors for \(r=5\) with varying \(s=6,\ldots,10\), \(\beta_\sigma=0, \dots, \sigma-1\), and \(n = \sigma^2 + \beta_\sigma\). Computed using 256-bit \textsc{BigFloat} and compared to 2048-bit \textsc{BigFloat}. Machine epsilon for 256-bit \textsc{BigFloat} is \(1.72\cdot10^{-77}\) (black dashed line). \textbf{Right:} Errors for \(r=25\) with varying \(s=25,\ldots,100\) and \(n=\sigma^2\), that is, \(\beta_\sigma=0\). Both the algorithm and a standard eigenvalue solver is used for \(T_{\sigma^2}(g_{25,s})\) for different precisions. Double precision (\textsc{Float64}) is used for computations with 53-bit precision, and \textsc{BigFloat} otherwise. Note the spikes down in the error when we have \(\gamma > 1\).}
    \label{fig:Error}
\end{figure}

Two sets of numerical experiments are performed to verify Conjecture~\ref{conj:1}. First, we have
\begin{align*}
    r&=1,\ldots,20,\\
    s&=r,\ldots,20,\\
    n&=\sigma^2+\beta_\sigma,\quad \sigma=r+s, \quad\beta_\sigma=0,\ldots, \sigma-1,
\end{align*}
that is, \(n\) grows as \(\sigma^2\), and the largest matrix is of size \(n=1639\).
We compute \(\lambda_j(T_n(g_{r,s}))\) using 256-bit \texttt{BigFloat} (machine epsilon \(\varepsilon=\mathcal{O}(10^{-77})\)) in \textsc{Julia} with the package \textsc{GenericSchur.jl}~\cite{genericschur.jl}. For select \(r\) and \(s\) were computed with 2048-bit precision for reference purposes.
With same precision we compute \(\lambda_+(T_n(g_{r,s}))\) using the proposed algorithm supporting Conjecture~\ref{conj:1}; algorithm implemented in Appendix~\ref{sec:appendix}. We present here a subset of these results in the left panel of Figure~\ref{fig:Error}, where we have plotted the error for \(r = 5\), \(s = \{6,\dots,10\}\). 

We define the error between the spectrum of $T_n(g_{r,s})$ computed using a standard solver and the presented algorithm as follows
\[ 
\max_{j=1,\ldots,n}\left|\left|\lambda_{\rho_1(j)}\left(T_n(g_{r,s})\right)_\text{std\vphantom{g}}\right|-\left|\lambda_{\rho_2(j)}\left(T_n(g_{r,s})\right)_\text{alg}\right|\right|, \]
where $\rho_1$ and $\rho_2$ are two permutations of $\{1,\ldots,n\}$ that respectively sort the eigenvalues of the computations in ascending order by absolute value. Furthermore, only the non-zero eigenvalue error is shown, as if we include the zero eigenvalues, then the error becomes the magnitude of the zero eigenvalues computed by the standard solver which is much larger than the error of the non-zero eigenvalues. 

We see that the general trend is for the error to decrease with \(\beta_\sigma\), that is that \(\beta_\sigma = 0\) will have the largest error; the condition number is also decreasing, exemplified earlier in Section~\ref{examples:2} and Figure~\ref{fig:charpoly12}. Thus, in the second experiment which goes over a larger range of \(\{r,s\}\) with a smaller \(n\), we will only compute for \(\beta_\sigma=0\) as this will yield the largest error. The second experiment performed was
\begin{align*}
    r&=1,\ldots,50,\\
    s&=r,\ldots,100,\\
    n&=3\sigma,\quad \sigma=r+s,
\end{align*}
where the largest matrix is \(n=450\). Computations are done with precisions \(\{53, 256, 512, 1024, 2048\}\) -bit for the full matrix \(T_n(g_{r,s})\) and the proposed algorithm where 53-bit is standard double precision (\textsc{Float64}) and the 2048-bit solution of the full matrix \(T_n(g_{r,s})\) is used as the reference solution. A subset of the results are shown in the right panel of Figure~\ref{fig:Error}. We show here how the error changes for a fixed value of \(r\) (\(r=25\) in this case), while ranging \(s=r,\ldots,100\). What we see is that the error increases as \(s\) increases for all precisions and this will be due to the ill-conditioned nature of the \(B_{n_\sigma}^{n,r,s}\) matrices. 

An additional numerical consideration is that the elements of \(B_{n_\sigma}^{n,r,s}\) can easily grow such that integer overflow is a concern; we use \textsc{Int128} as default in Appendix~\ref{sec:appendix} but \textsc{BigInt} may be used if this is an issue for large \(\{n,r,s\}\). Presented below is \(B_{3}^{231,38,39}\) where a standard 64-bit integer (\textsc{Int64}) would not work.
\begin{align*}
B_{3}^{231,38,39}=
\left[
\begin{array}{rrr}
2937189730080557577 & 9536995145808582886 & 11892438427558067162 \\
6599805415728025309 & 21429433573366650048 & 26722066585196691901 \\
5292633011830041853 & 17185071439388109015 & 21429433573366650048 \\
\end{array}
\right]
\end{align*}
The condition number of this matrix is \(4.8\cdot10^{46}\), highlighting the need to use high precision data types when computing the eigenvalues, \(\{1.0\cdot 10^{-27}, 2.2\cdot 10^7, 4.6\cdot 10^{19}\}\).

The final numerical aspect that we mention is the computational time of the algorithm. See Table~\ref{table:timings} which shows the minimal time used for the outlined algorithm in this paper and a standard eigenvalue solver (\texttt{eigvals} from \textsc{GenericSchur.jl} in \textsc{Julia}) for three different \(\{n,r,s\}\) triples. The timings were performed on a computer with an AMD Epyc 7502P CPU (32 core, 1024 GB RAM) and \textsc{BenchmarkTools.jl}~\cite{BenchmarkTools.jl-2016}.

The timings show, especially in the experiment with \(\{n,r,s\} = \{676, 7, 19\}\), that we can compute the eigenvalues of \(T_n(g_{r,s})\) with much higher precision using the presented algorithm with a comparable or faster execution time than a \textsc{Float64} computation of the full matrix (which will yield a considerable higher error). Optimizing the code in Appendix~\ref{sec:appendix}, and using iterative refinement~\cite{10.2307/2157168} we could probably increase the speed even further.
\begin{table}[H]
    \centering
    \caption{Timings, in seconds, using the algorithm of this paper and a standard eigenvalue solver for three \(\{n, r, s\}\) with different precisions. Standard double precision (\textsc{Float64}) is used for 53-bit and \textsc{BigFloat} for the other computations.}\label{table:timings}
    \begin{tabular}{r|rr|rr|rr}
    \toprule
    \multirow{2}{*}{precision (bit)} & \multicolumn{2}{c|}{\(\{n,r,s\} = \{49, 2, 5\}\)} & \multicolumn{2}{c|}{\(\{n,r,s\} = \{256, 4, 12\}\)} & \multicolumn{2}{c}{\(\{n,r,s\} = \{676, 7, 19\}\)} \\ \cline{2-7}
     & Algorithm & Standard & Algorithm & Standard & Algorithm & Standard \\ 
     \midrule
     53 & 0.00006 & 0.00055 & 0.00008 & 0.02654 & 0.00106 & 0.34084 \\
     64 & 0.00070 & 0.15268 & 0.00604 & 4.89205 & 0.01397 & 398.04027 \\
     128 & 0.00098 & 0.19567 & 0.00785 & 5.70696 & 0.01758 & 515.63287 \\
     256 & 0.00127 & 0.26923 & 0.01106 & 7.03421 & 0.02503 & 694.99100 \\
     512 & 0.00173 & 0.35154 & 0.01488 & 8.67083 & 0.03627 & 933.47496 \\
     1024 & 0.00284 & 0.60876 & 0.02403 & 14.36312 & 0.06258 & 1603.19599 \\
     2048 & 0.00577 & 1.23769 & 0.05347 & 23.06511 & 0.13496 & 2903.01721 \\ 
     \bottomrule
    \end{tabular}
\end{table}

\section{Conclusions}
\label{sec:conclusions}
Toeplitz matrices with two non-zero off-diagonals can be classified into three different cases. Firstly classical tridiagonal, and secondly ``symmetrically sparse tridiagonal''; we have analytical expressions for the eigenvalues (and eigenvectors) of these matrices. In the third, more general case we have conjectured in Conjecture~\ref{conj:1}, and argued in this paper, that it is possible to generate all the eigenvalues by generating the positive real subset of the eigenvalues of the original matrix from one or two smaller matrices (possibly with multiplicity), and then by rotation in the complex plane, and possibly adding zeros. 

Some aspects to further investigate, improve, and prove on the results of this paper are
\begin{itemize}
    \item Prove Conjecture~\ref{conj:1}.
    \item Devise an algorithm which works also for small matrices, see the restriction in subsection~\ref{sec:main:constructposreal}. Also, there might be some approach that does not yield any erroneous corners, given by the low rank perturbation \(R_{n_\sigma}^{n,r,s}\), at all.
    \item Devise an algorithm for computing the corresponding eigenvectors.
    \item The symbol \(\mathfrak{b}_{r,s}\) in \eqref{eq:shapirob} is not used in this article, since we are interested in the exact spectrum, and we do not know (in the general case) a grid \(\xi_j\), for \(j=1,\ldots, n_\sigma\), such that
    \( \lambda_j(B_{n_\sigma}^{n,r,s})=\mathfrak{b}_{r,s}(\xi_j)\), for \(j=1,\ldots, n_\sigma\); we do know this grid in Example 1 in Section~\ref{sec:examples}.
    However, we could use the symbol \(\mathfrak{b}_{r,s}(\theta)\) in conjunction with matrix-less methods~\cite{Barrera2018-tb,Bogoya2015-dt,Bogoya2017-ek,Bottcher2010-rf,Ekstr_m_2018,Ekstr_m_2017} to accurately and efficiently numerically approximate the spectrum of arbitrarily large $n_\sigma$. As an example, in~\cite{Ekstr_m_2021} we indeed have shown that the matrix-less method works on the matrices in Example 2, for \(\beta_\sigma=2\).
\end{itemize}

\section{Acknowledgements}
The first author would like to thank Albrecht Böttcher for his inspirational works in the field of Toeplitz-like matrices, and welcoming and illuminating discussions in Chemnitz.
We also thank Stefano Serra-Capizzano and Carlo Garoni for helpful suggestions for improving the paper.
The research of the second author is in part financed by the Centre for Interdisciplinary Mathematics (CIM) at Uppsala University.

\bibliographystyle{siam}
\bibliography{References}

\clearpage

\appendix

\section{Appendix: Algorithm}\label{sec:appendix}
In this Appendix we present a full implementation of the proposed algorithm in this paper, to compute the full spectrum of \(T_n(g_{r,s})\), assuming Conjecture~\ref{conj:1} is true. The algorithm is written with clarity in mind, and not optimised for efficiency.

\subsection{Setup and required packages}\label{code:packages}
\noindent
Required packages, and functions to construct square and rectangular Toeplitz matrices.

\noindent
\begin{minipage}{\linewidth}
\begin{lstlisting}[language=JuliaLocal, style=julia]
using LinearAlgebra
using GenericSchur # If using BigFloat (or other data types) instead of Float64
setprecision(BigFloat, 256) # If using BigFloat, choose a precision
# Toeplitz matrices of sizes n x n
function toeplitz(n::Integer, vc::Vector, vr::Vector, T::Type = Float64)
    Tn = zeros(T,n,n)
    for ii = 1:min(n,length(vc))
        Tn += vc[ii]*diagm(-ii+1=>ones(T,n-ii+1))
    end
    for jj = 2:min(n,length(vr))
        Tn += vr[jj]*diagm( jj-1=>ones(T,n-jj+1))
    end
    return Tn
end
# Rectangular Toeplitz matrices of sizes m x n
toeplitz(m::Integer, n::Integer, vc::Vector, vr::Vector, T::Type = Float64) = toeplitz(max(m,n), vc, vr, T)[1:m, 1:n]
\end{lstlisting}
\end{minipage}

\subsection{Construct \(B_{n_\sigma}^{n,r,s}\) matrices}\label{code:bns}
\begin{minipage}{\linewidth}
\begin{lstlisting}[language=JuliaLocal, style=julia]
function construct_B(n::Integer, r::Integer, s::Integer, T1::Type = Float64, T2::Type = Int128; remove_corner::Bool = true)
    σ = r+s
    β_σ = mod(n,σ)
    n_σ = div(n,σ)
    M, P = construct_M_P(r, s)
    B = I(n_σ)
    m = M[β_σ+1,:]
    p = P[β_σ+1,:]
    for kk = 1:r
        if p[kk] >= 0
            B *= transpose(construct_shifted_Tncm(n_σ,m[kk],T2)) * construct_shifted_Tncm(n_σ,1,T2)^p[kk]
        else
            B *= convert.(T2, transpose(construct_shifted_Tncm(n_σ,m[kk],T2)) / construct_shifted_Tncm(n_σ,1,T2))
        end
    end
    B = convert.(T1, B)
    # If beta>s and remove_corner is true, then extract corner from larger matrix and remove it from B
    if β_σ > s && remove_corner
        B_large = construct_B(σ^3 + β_σ, r, s, T1, T2, remove_corner = false) 
        n_min = min(r-1, n_σ)
        corner = abs.(B_large[1:n_min, end+1-n_min:end])
        B[1:n_min, end+1-n_min:end] += (-1)^isodd(n_σ)*corner
    end
    return B
end
\end{lstlisting}
\end{minipage}

\subsection{Construct \(\mathcal{M}\) and \(\mathcal{P}\)}\label{code:MP}

\begin{minipage}{\linewidth}
\begin{lstlisting}[language=JuliaLocal, style=julia]
function construct_M_P(r::Integer, s::Integer)
    σ = r+s
    m_array = repeat(1:σ, inner=r)
    p_array = repeat(-1:s, inner=r)
    M = toeplitz(σ, r, m_array[r:r+σ-1], ones(Int64,r), Int64)
    P = toeplitz(σ, r, reverse(p_array[2:σ+1]), p_array[σ+1:σ+r], Int64)
    # If r > 1, permute columns 
    if r != 1
        τ = mod(s, r)
        p_perm = mod.(τ:τ:r*τ, r)
        p_perm[end] = r
        M = M[:, p_perm]
        P = P[:, p_perm]
    end
    return M, P
end
\end{lstlisting}
\end{minipage}

\subsection{Construct \(T_n\left(e^{-\mathbf{i}\theta}c(\theta)^{m}\right)\) matrices}\label{code:tnc}
\begin{minipage}{\linewidth}
\begin{lstlisting}[language=JuliaLocal, style=julia]
function construct_shifted_Tncm(n::Integer, m::Integer, T::Type = Float64)
    v = binomial.(m,0:m)
    return toeplitz(n, v[2:end], reverse(v[1:2]), T)
end
\end{lstlisting}
\end{minipage}

\subsection{Compute positive real eigenvalues \(\lambda_+(T_{n}(g_{r,s}))\)}\label{code:realeigs}

\begin{minipage}{\linewidth}
\begin{lstlisting}[language=JuliaLocal, style=julia]
function construct_positive_real_eigenvalues(n::Integer, r::Integer, s::Integer, T1::Type = Float64, T2::Type = Int128)
    σ = r+s
    γ = gcd(r,s)
    σ_γ = div(σ, γ)
    β_γ = mod(n,γ)
    n_γ = div(n-β_γ,γ)
    r_γ = div(r,γ)
    s_γ = div(s,γ)
    # Trivial Case
    if n_γ < σ_γ
        n_γ_σ = div(n_γ,σ_γ)
        return repeat(zeros(T1, n_γ_σ), γ)
    end
    Bn_γ_σ = construct_B(n_γ,r_γ,s_γ,T1,T2)
    # Compute eigenvalues. If not using T1 = Float64, then it may be necessary to add a `maxiter=100000` argument to eigvals
    eBn_γ_σ = real.(Complex.(eigvals(Bn_γ_σ)).^(convert(T1, 1)/σ_γ))
    if β_γ == 0
        return sort(repeat(eBn_γ_σ, γ))
    else
        Bn_γp1_σ = construct_B(n_γ+1,r_γ,s_γ,T1,T2)
        eBn_γp1_σ = real.(Complex.(eigvals(Bn_γp1_σ)).^(convert(T1, 1)/σ_γ))
        return sort(vcat(repeat(eBn_γ_σ, γ-β_γ), repeat(eBn_γp1_σ, β_γ)))
    end
end
\end{lstlisting}
\end{minipage}

\subsection{Construct all eigenvalues \(\lambda_j(T_n(g_{r,s})\)) for \(j=1,\ldots,n\)}\label{code:alleigs}

\begin{minipage}{\linewidth}
\begin{lstlisting}[language=JuliaLocal, style=julia]
function construct_all_eigenvalues(n::Integer, r::Integer, s::Integer, T1::Type = Float64, T2::Type = Int128)
    σ = r+s
    γ = gcd(r,s)
    ω = div(σ,γ)
    positive_real_eigenvalues = construct_positive_real_eigenvalues(n,r,s,T1,T2)
    rotations = exp.((0:ω-1)*2im*convert(T1,pi)/ω)
    β_γ = mod(n, γ)
    n_γ = div(n-β_γ, γ)
    n_0 = (γ-β_γ)*mod(n_γ,ω) + β_γ*mod(n_γ+1,ω)
    return vcat(zeros(T1, n_0), kron(rotations, positive_real_eigenvalues))
end
\end{lstlisting}
\end{minipage}

\end{document}